\documentclass{conm-p-l}

\usepackage[cmtip,rotate,all]{xy}
\usepackage{url}
\newtheorem{theorem}{Theorem}[section]

\newtheorem{proposition}[theorem]{Proposition}

\theoremstyle{definition} \newtheorem{definition}[theorem]{Definition}
\newtheorem{example}[theorem]{Example}

\theoremstyle{remark} \newtheorem{remark}[theorem]{Remark}
\newtheorem*{further}{Further reading}

\numberwithin{equation}{section}

\makeatletter
\newcommand{\@cal}[1]{\expandafter\def\csname#1\endcsname{\mathcal{#1}}}
\newcommand{\@bb}[1]{\expandafter\def\csname
#1\endcsname{\mathbb{#1}}}
\newcommand{\@op}[1]{\expandafter\def\csname
#1\endcsname{\mathop{\operator@font {#1}}}}
\newcommand{\@rm}[1]{\expandafter\def\csname
#1\endcsname{\mathop{\operator@font {#1}}\nolimits}}
\newcommand{\romans}[1]{\@for\@i:=#1\do{\expandafter\@rm\expandafter{\@i}}}
\newcommand{\mathcals}[1]{\@for\@i:=#1\do{\expandafter\@cal\expandafter{\@i}}}
\newcommand{\mathops}[1]{\@for\@i:=#1\do{\expandafter\@op\expandafter{\@i}}}
\newcommand{\mathbbs}[1]{\@for\@i:=#1\do{\expandafter\@bb\expandafter{\@i}}}
\makeatother

\mathbbs{Z,N,L,P}
\newcommand{\CC}{\mathbb{C}}
\mathops{colim,holim,hocolim}
\mathcals{A,C,D,E,I,J,R,O,X,S,V,W}
\romans{ho,op,Map,Vect,sSet,Fun,Cones,Cat,Tor,Ext,%
Set,Hom,fr,StrCat,Top,sSet,id,Shv,Bord,Spec}

\newcommand{\oo}{\infty}
\newcommand{\io}{$(\oo,1)$}
\newcommand{\Catio}{\Cat_{(\oo,1)}}
\newcommand{\Catin}{\Cat_{(\oo,n)}}
\newcommand{\abs}[1]{\left|#1\right|}

\newcommand{\DELTA}{\mathbf{\Delta}}

\begin{document}

\title[Whirlwind Tour of $(\infty,1)$-categories]{A Whirlwind Tour of
the World of $(\infty,1)$-categories}

\author{Omar Antol\'{\i}n Camarena} \address{} \curraddr{}
\email{oantolin@math.harvard.edu} \thanks{}


\date{\today}

\begin{abstract}
  This introduction to higher category theory is intended to a give
  the reader an intuition for what \io-categories are, when they are
  an appropriate tool, how they fit into the landscape of higher
  category, how concepts from ordinary category theory generalize to
  this new setting, and what uses people have put the theory to. It is
  a rough guide to a vast terrain, focuses on ideas and motivation,
  omits almost all proofs and technical details, and provides many
  references.
\end{abstract}

\maketitle

\section{Introduction}

An \io-category is a category-like thing that besides objects and
morphisms has $2$-morphisms between morphisms, $3$-morphisms between
$2$-morphisms, and so on all the way to $\oo$; but in which all
$k$-morphisms for $k>1$ are ``invertible'', at least up to higher
invertible morphisms. This is the sort of invertibility that
homotopies have: the composition or concatenation of any homotopy with
its reverse is not actually the identity but it is homotopic to it. So
we can picture an \io-category as a ``homotopy theory'': a kind of
category with objects, morphisms, homotopies between morphisms, higher
homotopies between homotopies and so on.

Any context where there is a notion of homotopy, can benefit from the
use of \io-categories in place of ordinary categories. This includes
homotopy theory itself, of course, but also homological algebra and
more generally wherever Quillen's version of abstract homotopy theory,
the theory of model categories, has been used. Notions of homotopy are
perhaps more common than one might expect since the philosophy of
model categories shows that simply specifying a class of ``weak
equivalences'' in a category, a collection of morphisms which we wish
to treat as if they were isomorphisms, produces a notion of homotopy.
The theory of \io-categories plays a prominent role in derived
algebraic geometry, as can be expected from the very rough description
of the subject as being what is obtained by replacing the notion of
commutative rings in algebraic geometry by, say, commutative
differential graded algebras but only caring about them up to
quasi-isomorphism.

There are now several different formalizations or models of the notion
of \io-category, detailed comparison results relating the different
definitions and for one particular model of \io-category,
quasi-categories, a detailed account of how ordinary category theory
generalizes to the \io{} context \cite{Joyal}, \cite{HTT}, \cite{HA}.
(Many definitions and \emph{statements} of results from ordinary
category theory generalize straightforwardly to \io-categories, often
simply by replacing bijections of $\Hom$-sets with weak homotopy
equivalences of mapping spaces, but with current technology the
traditional \emph{proofs} do not generalize, and instead often require
delicate model specific arguments: most of this work has been done
using the model of quasi-categories.)

Giving a survey of the applications of (ordinary) category theory is
an impossible task: categories, and categorical constructions such as
products and adjoint functors, to give just two examples, appear in
very many fields of mathematics. Such a survey would turn into a
survey of much of mathematics. Writing an overview of the applications
of \io-categories could potentially be similarly doomed. This paper
attempts it anyway only because \io-categories are still relatively
new and have not fully caught on yet, making it possible to list a
reasonable portion of the current literature. Even so, this is just a
small entry point into the world of \io-categories and the broader
context of higher category theory.

The ideal reader of this survey is someone who has heard about
\io-categories (perhaps under the name $\oo$-categories), is
interested in reading some work that uses them (such as the derived
algebraic geometry literature), or is simply curious about them but
wishes to have a better idea of what they are and how they are used
before committing to read a rigorous treatment such as \cite{Joyal} or
\cite{HTT}. We will not assume any prior knowledge of \io-categories,
or even more than a cursory knowledge of 2-categories, but we will
assume the reader is comfortable with notions of ordinary category
theory such as limits, colimits, adjoint functors (but it's fine if the
reader can't give a precise statement of Freyd's Adjoint Functor
Theorem, for example). We also assume the reader is acquainted with
simplicial sets; if that's not the case we recommend reading
\cite{Friedman} as a gentle introduction that gives the basic
definitions and properties and focuses on conveying geometrical
intuition.\footnote{Anyone attempting to use the theory of
  \io-categories will need to know much more about simplicial sets,
  and would benefit from looking at a textbook such as \cite{MaySS} or
  \cite{GoerssJardine}.}

We will begin by briefly exploring the landscape of higher category
theory to give a context for \io-categories and describe some basic
guiding principles and requirements for the theory. Then we'll go on a
quick tour of all the different models available for \io-categories
and discuss the problem of comparing different definitions; an
exciting recent development is Barwick and Schommer-Pries's axiomatic
characterization of higher categories \cite{ClarkChris}. The next
section deals with practical aspects of working with \io-categories
and describes how concepts from ordinary category theory such as
isomorphisms, limits and colimits, adjunctions, monads, monoidal
categories and triangulated categories generalize to the \io{}
setting. The final section consists of (very!) brief descriptions
of some of the work that applies the theory of \io-categories.

\section{The idea of higher category theory}

The first hint of higher category theory comes from the category
$\Cat$ of categories. It not only has objects, which are categories,
and morphisms between them, functors, but there are also natural
transformations between functors. Indeed, $\Cat$ is the basic example
of a (strict\footnote{We'll explain what this means and have much more
to say about it in Section \ref{strict}.}) $2$-category, just as
$\Set$ is the basic example of a category. Of course, once we've
imagined, besides having objects and morphisms, having another layer
of things we'll call $2$-morphisms connecting the morphisms (in the
way natural transformations connect functors), there is no reason to
stop at $2$.

This gives us our first blurry picture of higher categories: an
$n$-category will have a collection of objects, and collections of
$k$-morphisms for $1 \le k \le n$ with specified identity morphisms
and composition operations for morphisms satisfying appropriate
associativity and unit axioms; an $\oo$-infinity category will be a
similar structure having $k$-morphisms for all $k \ge 1$.

\begin{remark} We are being very vague and purposefully so: there is a
large design space to explore. There are many possible forms for
composition laws and many ways of making the axioms precise, and there
are even many choices for the ``shape'' of morphisms, that is, choices
for what data specifies the analogue of domain-and-codomain of a
$k$-morphism. We won't have much to say about different shapes for
morphisms, so that discussion is postponed to section \ref{shapes}.
\end{remark}

Another way to visualize this idea is also already present in the
$2$-category $\Cat$: given two categories $\C$ and $\D$, $\Cat$
doesn't just have a \emph{set} of morphisms from $\C$ to $\D$, it has
a whole category $\Fun(\C, \D)$ whose objects are functors $\C \to \D$
and whose morphisms are natural transformations. Note the funny
re-indexing that takes place: \begin{enumerate} \item functors $\C \to
\D$ are $1$-morphisms in $\Cat$ but are $0$-morphisms (objects) in
$\Fun(\C,\D)$,

\item natural transformations are $2$-morphisms in $\Cat$ but are
$1$-morphisms in $\Fun(\C,\D)$. \end{enumerate}

This gives us an alternative inductive way to think of higher
categories: an $n$-category is like a $1$-category but instead of
having a $\Hom$-\emph{sets} between any pair of objects, its
$\Hom$-things are $(n-1)$-categories. Readers familiar with
\emph{enriched}\footnote{The definition of enriched category is
recalled in section \ref{strict}.} category theory will recognize that
this is similar to defining an $n$-category as a category enriched
over $(n-1)$-categories. That actually defines what is known as a
\emph{strict} $n$-category and we will argue in section \ref{strict}
that this notion does not capture the interesting examples that one
would want in higher category, so we really want some kind of category
``weakly enriched'' over $(n-1)$-categories. But before we discuss
that, notice that even brushing aside the issue of strictness, this
perspective does not help in defining $\oo$-categories, as the
inductive definition becomes circular in case $n = n-1 = \oo$.
However, if we restrict our attention to higher categories in which
above a certain level the morphisms behave like homotopies, we can
use the inductive perspective again.

Let's say an $(n,k)$-category is an $n$-category in which all
$j$-morphisms for $j\ge k+1$ are invertible in the sense homotopies
are: \emph{not} that every $j$-morphism $\alpha : x \to y$ has an
inverse $\beta$ for which $\beta \circ \alpha$ and $\alpha \circ
\beta$ are exactly \emph{equal} to the identity $(j+1)$-morphisms
$\id_x$ and $\id_y$, but only that there is a $\beta$ for which those
composites have invertible\footnote{Invertible in this same sense, so
this definition is recursive.} $(j+2)$-morphisms connecting them to
$\id_x$ and $\id_y$. Of course, if $j+2>n$, we \emph{do} require that
$\beta \circ \alpha = \id_x$ and $\alpha \circ \beta = \id_y$. In
other words, we can view any $n$-category as an $(n+1)$-category where
all the $(n+1)$-morphisms are identities. Finally, we can similarly
talk about $(\oo,k)$-categories (where reaching the top degree for
morphisms is not an issue), and the bulk of this survey will focus on
the \io{} case.

\begin{remark}\label{metaphor} A useful metaphor has us think of an
  invertible morphism between two objects as a proof that they are the
  ``same''. Just as with proofs of theorems in mathematics, sometimes
  one can argue that two proofs are ``really the same proof''; such an
  argument corresponds to an invertible $2$-morphism between two
  $1$-morphisms. Then we can think of proofs establishing that two
  ways of showing that two proofs that two objects are the same are
  the same are the same, and so on\footnote{Limited only by the number
    of times we are willing to say ``are the same'' in a row.}. In
  other words: an $(\oo,0)$-category, usually called an
  \emph{$\oo$-groupoid}, is what a set is forced to become if we are
  never satisfied to just note that two things can be proven to be the
  same, but instead we write down the proof and contemplate the
  possibility that different looking proofs can be proven to be the
  same. This is what people mean when they say higher category theory
  \emph{systematically replaces equality by isomorphism}.
\end{remark}

For these $(n,k)$-categories and $(\oo,n)$-categories, the inductive
perspective says that an $(n,k)$-category has $\Hom$-things which are
$(n-1,k-1)$-categories (which does not buy us anything new), but also
that an $(\oo,k)$-category has $(\oo,k-1)$-categories as
$\Hom$-things. To start picturing $(\oo,n)$-categories, we need to
know how to visualize $(\oo,0)$-categories, which is the next topic
on our agenda.

\subsection{The homotopy hypothesis and the problem with strictness}
\label{strict}

The $2$-category of categories is \emph{strict}, meaning that the
composition of both its $1$-morphisms and $2$ morphisms is associative
and has units (the identity $1$- and $2$-morphisms), as opposed to
being just something like ``associative up to homotopy''. One says
that in $\Cat$ composition is \emph{strictly} associative. As
mentioned above, it is easy to define $n$-categories with strictly
associative and unital compositions inductively using the notion of
enriched category. Recall that given a monoidal category $\V$ with
tensor product given by a functor $\otimes : \V \times \V \to \V$, a
category $\C$ enriched over $\V$ (sometimes called a $\V$-category)
consists of \begin{itemize} \item a collection of objects, \item
$\Hom$-objects $\C(X,Y) \in \V$, for every pair of objects of $\C$,
\item composition morphisms $\C(Y,Z) \otimes \C(X,Y) \to \C(X,Z)$ of
$\V$, for every triple of objects of $\C$, \item identities given as
morphisms $I \to \C(X,X)$ in $\V$, for every object of $\C$ (where $I$
is the tensor unit in $\V$), \end{itemize} and this data is required
to satisfy obvious unit and associativity axioms (whose precise
statement requires using the unit and associativity constraints of
$\V$). When $\V$ is a category with finite products, we can take the
tensor product to be the categorical product (and $I$ to be the
terminal object); when equipped with this tensor product, $\V$ is said
to be a \emph{Cartesian} monoidal category. There is also a notion of
$\V$-enriched functor between two categories $\C$ and $\D$ enriched
over $\V$: a function associating to every object $X \in \C$ an object
$FX \in \D$, plus a collection of morphisms of $\V$, $\C(X,Y) \to
\D(FX,FY)$ compatible with identities and composition.

We can now give the inductive definition of strict $n$-categories:

\begin{definition} A \emph{strict $n$-category} is a category enriched
over the Cartesian monoidal category $\StrCat_{n-1}$. The category
$\StrCat_n$ whose objects are all strict $n$-categories and whose
morphisms are $\StrCat_{n-1}$-enriched functors is easily seen to have
finite products, making the recursion well defined. The base case can
be taken to be $\StrCat_1$, the ($1$-)category of categories and
functors or even $\StrCat_0 = \Set$.\footnote{In fact, one can make
sense of $\StrCat_{-1}$ and $\StrCat_{-2}$ as well! It's left as a fun
exercise for the reader.} \end{definition}

The only higher category we've mentioned so far is $\Cat$, and it
\emph{is} a strict $2$-category, but that's more or less it for
naturally occurring examples of strict $2$-categories, in the sense
that almost all natural examples have an air about them of functions
and composition of functions.

\begin{example} A monoid $M$ can be regarded as a category that has a
single object $x$ for which $\Hom(x,x) = M$ with composition given by
the monoid multiplication. In a similar way one cat try to turn a
monoidal category $\V$ into a $2$-category with one object $x$ for
which $\Hom(x,x) = \V$ with composition given by the tensor product in
$\V$. This does not produce a strict $2$-category unless the tensor
product is strictly associative and unital. The point of this example
is that most naturally occurring monoidal categories are \emph{not}
strict. For example, the tensor product of vector spaces is only
defined up to canonical isomorphism, and while $(U \otimes V) \otimes
W \cong U \otimes (V \otimes W)$, it is exceedingly unlikely that any
actual choice of specific vector spaces for all tensor products would
render both sides exactly \emph{equal}. Similar remarks apply to
products for Cartesian monoidal categories. \end{example}

\begin{remark} There is a standard notion of non-strict $2$-category: the
notion of \emph{bicategory} due to B\'{e}nabou \cite{Benabou} (or see
\cite{Lack}), that has a definition very similar to the usual
definition of monoidal category and which reduces to it in the case of
a bicategory with a single object. \end{remark}

While we have given what we feel are natural examples of
$2$-categories that fail to be strict, maybe they do not make a
conclusive case for the need to weaken the associativity and unitality
axioms: MacLane's coherence theorem for monoidal categories shows
that any monoidal category is (monoidally) equivalent to one where the
tensor product is strictly associative. And more generally
any $2$-category is equivalent to a strict one\footnote{On the
other hand maybe we do have a conclusive case for considering more
general notions than strict \emph{functors}: not every functor of
bicategories between strict $2$-categories is equivalent to a strict
$2$-functor! See \cite[Lemma 2]{LackTri} for an example. B\'{e}nabou
has expressed the view that the point of bicategories is not that they
are non-strict themselves, but that they are the natural home for
non-strict functors.} (see \cite{LeinsterBicat} for an expository account).
But once we get to $3$-categories, the situation is different: there are
examples that cannot be made strict. We'll give an explicit example in
section \ref{sphere}, namely, the fundamental $3$-groupoid of $S^2$;
but first we will discuss fundamental higher groupoids and their role
in higher category theory.

Higher groupoids are special cases of higher categories, namely an
$n$-groupoid, in the terminology explained above, is an
$(n,0)$-category and an $\oo$-groupoid is a $(\oo,0)$-category. Before
we explain what higher fundamental groupoids should be, recall that
the fundamental groupoid packages the fundamental groups of a space
$X$ at all base points into a single category $\pi_{\le 1} X$ whose
objects are the points of $X$ and whose morphisms $x \to y$ are
endpoint-preserving homotopy classes of paths from $x$ to $y$.
Composition is given by concatenation of paths (which is not strictly
associative and unital before we quotient by homotopy\footnote{We could
use Moore paths, which are maps $[0,\ell] \to X$ for some $\ell \ge 0$
called the length of the path. When concatenating Moore paths, the lengths
add. This operation \emph{is} strictly associative and unital, but (1) the
category of Moore paths is not a groupoid, since the reversal of a path
only is an inverse up to homotopy, and (2) there is no analogue of
Moore paths for the fundamental $n$-groupoid when $n>1$.}). For a space
$X$ that has some non-zero higher homotopy groups, $\pi_{\le 1}X$
clearly does not contain all the homotopical information of $X$, but
for \emph{$1$-types} it does.

\begin{definition} A space\footnote{For technical reasons, space here
shall mean ``space with the homotopy type of a CW-complex'', otherwise
some of the statements need homotopy equivalences replaced by weak
homotopy equivalence.} $X$ is called an $n$-type if $\pi_k(X,x) = 0$
for all $k > n$ and all $x \in X$. \end{definition}

The homotopy theory of $1$-types is completely captured by groupoids:

\begin{enumerate}
  \item The fundamental groupoid functor induces an equivalence
    between (a) the homotopy category of $1$-types, where the
    morphisms are homotopy classes of continuous functions between
    $1$-types\footnote{This definition of the category is correct
      because we took $1$-types to have the homotopy type of a
      CW-complex; we could instead consider the category obtained from
      $1$-types by inverting weak homotopy equivalences.}, and (b) the
    homotopy category of groupoids, whose morphisms are equivalence
    classes of functors between groupoids, two functors being
    equivalent if there is a natural isomorphism between them.

  \item The inverse of the equivalence described above can be given by
    a \emph{classifying space} functor $B$ that generalizes the
    well-known construction for groups and is defined before passing
    to homotopy categories, i.e., is a functor from the category of
    groupoids to the category of $1$-types. Any groupoid $G$ is
    equivalent to $\pi_{\le 1}BG$, and any $1$-type $X$ is homotopy
    equivalent to $B \pi_{\le 1} X$.

  \item Given two $1$-types $X$ and $Y$ (or more generally, an
    arbitrary space $X$ and a $1$-type $Y$), the space of maps
    $\Map(X,Y)$ is a $1$-type and its fundamental groupoid is the
    category of functors $\Fun(\pi_{\le 1} X, \pi_{\le 1} Y)$ (which
    is automatically a groupoid too).
\end{enumerate}

This means that homotopy theoretic questions about $1$-types can be
translated to questions about groupoids which thus provide complete
algebraic models for $1$-types. This is the simplest case of perhaps
the main guiding principle in the search for adequate definitions
in higher category theory: the \emph{homotopy hypothesis} proposed by
Grothendieck in \cite{PursuingStacks}. As is common now, we interpret
(and phrase!) it as stating desired properties of a theory of higher
categories.

\medskip

\textbf{The homotopy hypothesis:} Any topological space should have a
fundamental \emph{$n$-groupoid} for each $n$ (including $n = \oo$).
These should furnish all examples of $n$-groupoids in the sense that
every $n$-groupoid should be equivalent to the fundamental
$n$-groupoid of some space. Furthermore, the theory of $n$-groupoids
should be the ``same'' as the homotopy theory of $n$-types (where if
$n=\oo$, ``the homotopy theory of $n$-types'' is just ``homotopy
theory'').

\medskip

Notice that this only puts requirements on $(n,k)$-categories for
$k=0$, so it certainly does not tell the whole story of higher
category theory, but it is enough to rule out basing the theory on
strict $n$-categories as we'll see in the next section. This means
that we must search for definitions of higher categories that are
\emph{non-strict} or \emph{weak}, in the sense we mentioned
monoidal categories are weak: instead of associativity meaning that
given three $k$-morphisms $f$, $g$ and $h$, the composites
$(f \circ g) \circ h$ and $f\circ (g \circ h)$ are
\emph{equal}, we should only require them to be linked by an
invertible $(k+1)$-morphism $(f \circ g) \circ h \to f\circ (g \circ h)$
that could be called an \emph{associator}.
The reader familiar with the definition of monoidal category will
know that these associators should satisfy a condition of their own.
Given four $k$ morphisms $f$, $g$, $h$ and $k$, we can relate the
composites $((fg)h)k$ and $f(g(hk))$ in two different ways (we've
dropped the $\circ$ for brevity):

\centerline{
  \xymatrix @C=-3mm{
     && (f(gh))k \ar[drr] && \\
     ((fg)h)k \ar[urr] \ar[dr] &&&& f((gh)k) \ar[dl] \\
     & (fg)(hk) \ar[rr] && f(g(hk)) & \\
  }
}

For the case of monoidal categories (where $f$, $g$, $h$ and $k$ are
objects, $\circ = \otimes$, and the associator is a $1$-morphism)
we've reached the top level already and we require this diagram to
commute; but in a higher category we can instead requires this to
commute only up to an invertible $(k+2)$-morphism we could call a
\emph{pentagonator}.  This pentagonator must satisfy its own
condition, but only up to a higher morphisms and so on. This kind of
data ---the associators, pentagonators, etc.--- are what is meant to
exist when saying an operation is associative up to \emph{coherent
  homotopy}.

Clearly, drawing these diagrams gets complicated very quickly
and indeed, definitions of $n$-categories along these lines have
only been written down for $n$ up to $4$ ---for a definition of
tricategories see \cite{GPS} or \cite{Gurski}, for tetracategories
see \cite{Trimble} or \cite{Hoffnung}. Instead people find ways
of implicitly providing all these higher homotopies in a clever
roundabout way. We'll see some examples in the section on models
of \io-categories.

\subsection{The $3$-type of $S^2$}\label{sphere}

We will show that the fundamental $3$-groupoid of $S^2$ is not
equivalent to a strict $3$-groupoid, or, in other words, that there is
no strict $3$-groupoid that models the $3$-type of $S^2$, which is
commonly denoted $P_3 S^2$ in the theory of Postnikov
towers.\footnote{Recall that the $3$-type of $S^2$ can be obtained,
say, by building $S^2$ as a CW-complex and then inductively attaching
cells of larger and larger dimension to kill all homotopy groups
$\pi_i$ for $i \ge 4$.} What we mean by ``models'' is that we assume
the existence of classifying space functors (with certain properties
we'll spell out later that \emph{are} satisfied for the ``standard
realization functors'', see the discussion after Theorem 2.4.2 of
\cite{Simpson}) that produce an $n$-type $BG$ for a strict
$n$-groupoid, and we say $G$ models a space $X$ if $BG$ is homotopy
equivalent to $X$. The argument shows, more generally, that if $X$ is
a simply connected $n$-type modeled by a strict $n$-groupoid $G$, $X$
is in fact an infinite loop space and even a product of
Eilenberg-MacLane spaces.

Let's investigate when we can \emph{deloop} a given strict
$n$-groupoid $G$, i.e., when $G$ can be realized as $\Hom_H(x,x)$ for
some strict $(n+1)$-groupoid $H$ with a single object $x$. That's easy
enough: if there exists such an $H$, $G$ inherits from it a
composition $\mu : G \times G \to G$ which makes it into a monoid
object in $\StrCat_n$, and clearly for each such monoid structure we
can form a delooping $H$. If we want to deloop more than once, we need
a monoid structure on $H$. And here something remarkable happens: a
monoid structure $\nu : H \times H \to H$, in particular restricts to
a new monoid structure $\nu_G : G \times G \to G$ on $G =
\Hom_H(x,x)$, and, since $\nu$ is a $\StrCat_n$-enriched functor, this
$\nu_G$ must be compatible with composition in $H$, that is, with
$\mu$. The end result is that $G$ has \emph{two} monoid structures one
of which is homomorphism for the other. The classical Eckmann-Hilton
argument\footnote{This says that if a set $M$ has two different monoid
structures given by products $\cdot$ and $\ast$, and we have $(a \cdot
b) \ast (c \cdot d) = (a \ast c) \cdot (b \ast d)$ ---which says that
$\ast : (M,\cdot) \times (M,\cdot) \to (M,\cdot)$ is a monoid
homomorphism--- then $\cdot =\ast$ and $M$ is commutative. Here we are
actually using an extension to strict $n$-categories instead of sets,
which is essentially obtained by applying the classical statement to
each degree of morphism separately.} implies that $\nu_G = \mu$ and
that they are commutative. Also, conversely, if $\mu$ is a commutative
monoid structure for $G$, $H$ does in fact become a commutative monoid
under $\nu(x) = x$, $\nu_G = \mu$.

This means that in the world of strict $n$-groupoids, delooping twice
is already evidence that you can deloop arbitrarily many times! Using
this it is easy to see why you can't find a strict $3$-groupoid that
models $P_3 S^2$. If there were such a groupoid $G$, without loss of
generality we could assume $G$ had a single object $x$ and a single
$1$-morphism $\id_x$: otherwise just take the sub-strict-$3$-groupoid
consisting of $x$, $\id_x$ and the groupoid $\Hom_G(\id_x, \id_x)$.
But then $G$ is the second delooping of $\Hom_G(\id_x, \id_x)$, which
shows that $G$ in turn can be delooped arbitrarily many times. If we
had classifying spaces for groupoids that were compatible with looping
(by which we mean we had an $n$-type $BG$ for each strict $n$-groupoid
$G$ such that if $G$ has a single object $x$, $\Omega BG$ is weakly
homotopy equivalent to $B(\Hom_G(x,x))$), it would follow that $P_3
S^2$ is an infinite loop space, which it is not. In fact, if
classifying spaces preserved products (i.e., $B(G \times H) \cong BG
\times BH$), we'd have that $P_3 S^2$ would be a topological abelian
monoid and thus homotopy equivalent to a product of Eilenberg-MacLane
spaces. It would then have to be $K(\Z,2) \times K(\Z,3)$, but it is
not, since, for example, the Whitehead product $\pi_2 S^2 \times \pi_2
S^2 \to \pi_3 S^2$ is non-zero.

\begin{remark} Vanishing of the Whitehead product $\pi_2 \times \pi_2
\to \pi_3$ does not guarantee that a $3$-type can be modeled by a
strict $3$-groupoid. Consider the space $X = P_3 Q S^2 = P_3 \colim
\Omega^n \Sigma^n S^2$ whose Whitehead product is $0$ simply for
torsion reasons: $\pi_2 P_3QS^2 = \Z/2$, the first stable homotopy
group of spheres. One can see $X$ is not homotopy equivalent to
$K(\Z,2) \times K(\Z,3)$ by looking at the operation $\pi_2 W \to
\pi_3 W$ given by composing (maps representing homotopy classes) with
the generator of $\pi_3(S^2)$: this operation is non-zero for $W=X$,
but is zero for a product of Eilenberg-MacLane spaces. By the argument
above, $X$ is not modeled by a strict $3$-groupoid. \end{remark}

\begin{further} Carlos Simpson \cite[Section 2.7]{Simpson} proved that
there is no classifying space functor for strict $3$-groupoids such
that $BG$ is homotopy equivalent to $P_3S^2$ under weaker assumptions
than we sketched above: he does not assume that classifying spaces are
compatible with looping, in fact, he does not require there to be a
family of classifying space functors for strict $n$-groupoids for all
$n$ at all; just a single functor for $n=3$ satisfying the minimal
requirements that $BG$ be a $3$-type and that the homotopy groups of
$BG$ are functorially isomorphic to algebraically defined ones for
$G$. The simpler argument we sketched (under the stronger assumption
of compatibility with looping) can be found in \cite[Section
2.6]{Simpson}. Clemens Berger proved a stronger result characterizing
all connected $3$-types (not necessarily simply connected) that can be
modeled by strict $3$-groupoids \cite[Corollary 3.4]{Berger}.
\end{further}

\subsection{Other shapes for cells}\label{shapes}

There are other possibilities for the shapes of the morphisms in an
$n$-category, of which we'll give a brief representative list. Here we
will call the morphisms \emph{cells}, since the word ``morphism'' is a
little awkward when a $k$-morphism does not simply go from one
$(k-1)$-morphism to another.

In the case of the $2$-category of categories, the $2$-morphisms,
which are natural transformations, go between two $1$-morphisms
(functors) that are parallel, i.e., that share their domain and share
their codomain. This pattern can be generalized for higher morphisms
and is called \emph{globular}, because drawings of such morphisms
looks like topological balls, or more precisely like one of the usual
CW-complex structures on disks: the one in which the boundary of the
disk is divided into hemispheres meeting along a sphere, which is also
divided into hemispheres and so on. Here's a picture of a globular
$2$-cell:

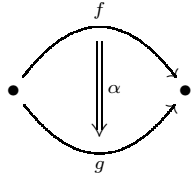
\centerline{ \xymatrix{ \bullet \ar@/^2pc/[rr]_{\quad}^{f}="1"
\ar@/_2pc/[rr]_{g}="2" && \bullet \ar@{}"1";"2"|(.2){\,}="7"
\ar@{}"1";"2"|(.8){\,}="8" \ar@{=>}"7" ;"8"^{\alpha} }}

For another example of a shape for morphisms think of homotopies,
homotopies between homotopies, and so on. As we mentioned in the
introduction, this is one of the examples we are trying to capture.
Such higher homotopies are maps $X \times [0,1]^n \to Y$, and so
naturally have a \emph{cubical} shape. A cubical $2$-morphism looks
like a square, and its analogue of domain-and-codomain, the boundary
of the square, has four objects and four $1$-morphisms:

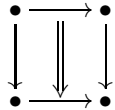
\centerline{ \xymatrix @C=3mm { \bullet \ar[d] \ar[rr] & \ar@{=>}[d] &
\bullet \ar[d] \\ \bullet \ar[rr] & & \bullet } }

When we get to discussing models for \io-categories, and specifically
the model of quasi-categories (which are simplicial sets satisfying
some condition), we will encounter another shape for morphisms:
\emph{simplicial}. A $2$-morphism is shaped like a triangle:

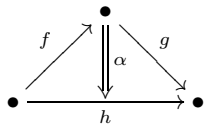
\centerline{ \xymatrix{ & \bullet \ar[dr]^{g} \ar@{=>}[d]^-{\alpha} &
\\ \bullet \ar[rr]_{h} \ar[ur]^{f} & & \bullet}}

We can interpret $\alpha$ as a homotopy between $g\circ f$ and $h$ (or
alternatively we can interpret composition as multivalued, in which
case, $h$ is some composite of $g$ and $f$, and $\alpha$ is a witness
to that fact). Similarly we can think of a higher dimensional simplex
as being a coherent collection of homotopies between composites
of a string of $1$-morphisms. See section \ref{qcat}.

There are more elaborate cell shapes as well, such as \emph{opetopes},
introduced by John Baez and James Dolan \cite{BaezDolan3} (see also
\cite{BaezIntro} which besides describing opetopes and the proposed
definition of higher category based on them, is a nice introduction to
$n$-categories generally). These can be interpreted as being a
homotopy between the result of evaluating a \emph{pasting diagram} and
a specified target morphism. This is analogous to the above
interpretation of simplices, but allowing for more general pasting
diagrams than those given by strings of $1$-morphisms.

\begin{further}
  To look at pictures of the zoo of higher categories, we recommend
  the illustrated guidebook by Eugenia Cheng and Aaron Lauda
  \cite{ChengLauda}. (In particular, the above description of opetopes
  is meaningless without pictures, which can be seen there.)  For a
  concise list of many of the available definitions for $n$-category
  and $\oo$-category see \cite{LeinsterDef}. See also the book
  \cite{LeinsterBook}, particularly Chapter 10.

  Those sources concentrate on definitions attempting to capture
  $n$-categories without any requirement of invertibility of morphisms.
  Thanks to the homotopy hypothesis and the availability of topological
  spaces, simplicial sets and homotopy theory it has turned out somewhat
  easier in practice to work with notions of $(\oo,n)$-categories
  (which of course include $n$-categories as a special case). As John
  Baez said about climbing up the categorical ladder from $1$-groupoids
  to $\oo$-groupoids \cite{TWF223}\footnote{John Baez's web column
  \emph{This Week's Finds} is a highly recommended source for intuition
  about higher categories.}:

  \begin{quotation}
    [\dots] the $n$-category theorists meet up with the topologists ---
    and find that the topologists have already done everything there is
    to do with $\oo$-groupoids\dots{} but usually by thinking of them
    as spaces, rather than $\oo$-groupoids!

    It's sort of like climbing a mountain, surmounting steep cliffs with
    the help of ropes and other equipment, and then finding a Holiday
    Inn on top and realizing there was a 4-lane highway going up the
    other side.
  \end{quotation}

  For the homotopical perspective and a focus on $(\oo,n)$-categories
  see \cite{Simpson}. The rest of this survey will mostly focus on
  \io-categories.
\end{further}

\subsection{What does (higher) category theory do for us?}

The reader might be asking now how exactly higher category theory is
useful in mathematics. Here is one possible answer, a purely
subjective and personal answer, and should be disregarded if the
reader does not find it convincing. It is now widely recognized that
category theory is a highly versatile and profitable organizing
language for mathematics. Many fields of mathematics have objects of
interest and distinguished maps between them that form categories,
many comparison procedures between different kinds of objects can be
represented as functors and, perhaps, most importantly, basic notions
from category theory such as products, coproducts (or more general
limits and colimits) and adjoint functors turn out to be well-known
important constructions in the specific categories studied in many
fields. While it is not reasonable to expect that category theory will
swoop in and solve problems from other fields of mathematics, phrasing
things categorically does help spot analogies between different fields
and to pinpoint where the hard work needs to happen: often arguments
are a mix of ``formal'' parts, which depend very little on the
detailed structure of the objects being studied, and ``specific''
parts which involve understanding their distinguishing properties;
categorical language makes short work of many formal arguments, thus
highlighting the remainder, the ``essential mathematical content'' of
an argument. Higher category theory promises to extend the scope of
such formal methods to encompass situations where we wish to consider
objects up to a weaker notion of equivalence than isomorphism; for
example, we almost always wish to consider categories up to
equivalence, in homological algebra we consider chain complexes up to
quasi-isomorphism, and in homotopy theory we consider space up to
homotopy equivalence or weak homotopy equivalence.

\section{Models of \io-categories}\label{io-models}

This sections gives a list of the main models of \io-categories and
attempts to motivate each definition. We spend more time discussing
quasi-categories than the other models, and in later sections we'll
mostly just use quasi-categories whenever we need particular models.
The reader will notice an abundance of simplicial sets appearing in
the definitions, and is warned again that some basic knowledge of them
will be required.

Ideally we would describe for each model, say, \begin{itemize} \item
the definition of \io-category, \item the corresponding notion of
functor and even the \io-category $\Fun(\C,\D)$ of functors between
two given \io-categories, \item how to retrieve the
$\Hom$-$\oo$-groupoid, or \emph{mapping space} $\Map_\C(X,Y)$ between
two objects of a given \io-category, \item the homotopy category $\ho
\C$ of a given \io-category, which is the ordinary category with the
same objects as $\C$ and whose morphisms correspond to
homotopy\footnote{Recall that since in an \io-category $2$-morphisms
and higher are invertible, we tend to think of them as homotopies}
classes of morphisms in $\C$. \end{itemize}

Sadly, for reasons of space we will not do all of those for each
model, but we hope to mention enough of these to give an idea of how
the story goes.

One excellent feature of the \io{} portion of higher category theory
is that the problem of relating different definitions has a
satisfactory answer which will be described in the following section.

\begin{further} For a more detailed introduction to the different
models and the comparison problem, we recommend \cite{BergnerSurvey},
\cite{JoyalTierney} or \cite{Porter}. \end{further}

\subsection{Topological or simplicial categories}

As we mentioned above, we can think of an \io-category as a category
weakly enriched in $\oo$-groupoids, and to satisfy the homotopy
hypothesis we could ``cheat'' and define $\oo$-groupoids as
topological spaces or simplicial sets (whose homotopy theory is
well-known to be equivalent to that of topological spaces). It turns
out that one can always ``strictify'' the enrichment in
$\oo$-groupoids, meaning that we can model \io-categories using:

\begin{definition} A \emph{topological category} is a category
enriched over the category of topological spaces\footnote{Instead of
the category of all topological spaces it is better to use a so-called
``convenient category of spaces'' \cite{Steenrod}, such as compactly
generated weakly Hausdorff spaces (see, for instance,
\cite{Strickland}). This is to make the comparison with other models
smoother and is a technical point the reader can safely ignore.}. A
\emph{simplicial category} is similarly a category enriched over the
category of simplicial sets. \end{definition}

These models of \io-categories are perhaps the easiest to visualize
and are a great psychological aid but are inconvenient to work with in
practice because, among other problems, enriched functors do not
furnish all homotopy classes of functors between the \io-categories
being modeled, unless the domain and codomain satisfy appropriate
conditions\footnote{Namely, that the domain be cofibrant and the
codomain be fibrant in the model structures discussed in section
\ref{comparison}.}.

Finally, notice that although we ``cheated'' by putting the homotopy
hypothesis into the definition, there is a sense in which we don't
trivially get it back out! We obtained a definition of \io-category
through enrichment from a definition of $\oo$-groupoid, but having
done so we now have a second definition of $\oo$-groupoid: an
\io-category in which all $1$-morphisms are invertible (up to higher
morphisms, as always). In terms of the homotopy category $\ho\C$, this
definition of $\oo$-groupoid that accompanies any notion of
\io-category is simply: an \io-category $\C$ for which $\ho\C$ is a
groupoid.

For topological or simplicial categories it is easy to construct $\ho
\C$: take the set of morphisms between $X$ and $Y$ in $\ho\C$ to be
$\pi_0(\C(X,Y))$; since $\pi_0$ preserves products, the composition
law in $\C$ induces a composition for $\ho\C$. Now, given a
topological\footnote{The simplicial case is analogous.} category $\C$
for which $\ho\C$ is a groupoid, what is the space $C$ that this
$\oo$-groupoid is supposed to correspond to? Think first of the case
when $\C$ has a single object $X$. Then $M:=\C(X,X)$ is a topological
monoid and $\ho\C$ being a groupoid just says that $\pi_0(M)$ is a
group under the operation induced by the multiplication in $M$. The
topological category $\C$ is a delooping\footnote{See section \ref{sphere}.}
 of $M$, so we should have
$\Omega C \cong M$, and there is such a space: the classifying
space $C=BM$ of $M$; when $\pi_0(M)$ is a group, the unit map
$M \to \Omega BM$ is weak homotopy equivalence. For general
groupoids $\ho\C$, the space $C$ corresponding to $\C$ will be
a disjoint union of classifying spaces of the monoid of endomorphisms
of an object chosen from each component of $\ho\C$.

\subsection{Quasi-categories}\label{qcat}

There are two classes of examples we certainly wish to have in any
theory of \io-categories: (a) ordinary categories (just add identity
morphisms in all higher degrees), and (b) $\oo$-groupoids, which by
the homotopy hypothesis we can take to be anything modeling all
homotopy types of spaces. After spaces themselves, the best known
models for homotopy types are \emph{Kan complexes}, simplicial sets
$X$ that satisfy the \emph{horn filler} condition: that every map
$\Lambda^n_k \to X$ extends to a map $\Delta^n \to X$. (Recall that
$\Lambda^n_k$ is obtained from the boundary $\partial \Delta^n$ of
$\Delta^n$ by removing the $k$-th face.) Also, every category $\C$ has
a nerve which is a simplicial set whose $n$-simplices are indexed by
strings of $n$ composable morphisms of $\C$; and the nerve functor $N
: \Cat \to \sSet$ is fully faithful. So inside the category $\sSet$ of
simplicial sets we find both ordinary categories and Kan complexes and
so we might expect to find a good definition of an \io-category as a
special kind of simplicial set. The following easy characterization of
those simplicial sets which arise as nerves of categories shows what
to do:

\begin{proposition}\label{nerves}
  A simplicial set $X$ is isomorphic to the nerve of
  some category if and only if every map $\Lambda^n_k \to X$ with
  $0<k<n$ extends \emph{uniquely} to a map $\Delta^n\to X$.
\end{proposition}

The least common generalization of the condition above and the
definition of Kan complex is:

\begin{definition} A \emph{quasi-category} is a simplicial set in
which all \emph{inner horns} can be filled, that is, in which every
map $\Lambda^n_k \to X$ with $0<k<n$ extends to a map $\Delta^n \to
X$. \end{definition}

Probably the greatest advantage of quasi-categories over other models
for \io-categories is how straightforward it is to deal with functors.
A functor $\C \to \D$ between two quasi-categories is simply a map of
simplicial sets: the structure of the quasi-categories makes any such
maps behave like a functor. (This is related to the nerve functor
being fully faithful.) Moreover, there is a simple way to obtain the
\io-category of functors between two quasi-categories: it is just the
simplicial mapping space\footnote{This is the internal $\hom$ in
$\sSet$, its $n$-simplices are simplicial maps $\C \times \Delta^n \to
\D$.} $\D^\C$, which is automatically a quasi-category whenever $\C$
and $\D$ are. In fact, more generally, given a quasi-category $\C$,
and an arbitrary simplicial set $X$, $\C^X$ is a quasi-category which
we think of as the category of $X$-shaped diagrams in $\C$.

The definition of quasi-category is very clean, but it may seem
mysterious that it does not mention anything like composition of
morphisms. Quasi-categories have something like a ``multivalued''
composition operation. Consider two morphisms $f : X \to Y$ and $g : Y
\to Z$ in a quasi-category $\C$ ---this really means that $X$, $Y$ and
$Z$, are vertices or $0$-simplices in the simplicial set $\C$ and that
$f$ and $g$ are $1$-simplices with the specified endpoints. The data
$(X,f,Y,g,Z)$ determines a map $\Lambda^2_1 \to \C$, that we display
by drawing $\Lambda^2_1$ and labeling the simplices by their images in
$\C$. A filler for this horn is a $2$-simplex $\alpha$ whose third
edge $h$ gives a possible composite of $g$ and $f$. The $2$-simplex
itself can be considered to be some sort of certificate that $h$ is a
composite of $g$ and $f$. There may be more than one composite $h$,
and for a given $h$ there may be more than one certificate.

\centerline{$ \begin{gathered} \xymatrix @R=6mm @C=2mm { & Y
\ar[dr]^{g} & \\ X \ar[ur]^{f} & & Z } \end{gathered}
\hspace{7mm}\leadsto\hspace{7mm} \begin{gathered} \xymatrix @R=6mm
@C=2mm { & Y \ar[dr]^{g} \ar@{=>}[d]_-{\alpha}& \\ X \ar[ur]^{f}
\ar[rr]_{h} & & Z } \end{gathered}$ }

This might seem like chaos, but homotopically composition is
well-defined in a sense we'll now make precise. The space of
composable pairs of $1$-simplices in $\C$ is given by the simplicial
mapping space $\C^{\Lambda^2_1}$ and the space of ``certified
compositions'' is similarly $\C^{\Delta^2}$. The set of vertices of
$\C^{\Lambda^2_1}$ is precisely the set of pairs of composable
$1$-simplices, and the higher dimensional simplices capture homotopies
between diagrams of composable pairs, and homotopies between those,
and so on. Similar remarks apply to $\C^{\Delta^n}$.

\begin{proposition}[Joyal]\label{comp}
  For a quasi-category $\C$, the map $\C^{\Delta^2} \to
  \C^{\Lambda^2_1}$ induced by composition with the inclusion
  $\Lambda^2_1 \hookrightarrow \Delta^2$ is a trivial Kan fibration,
  which implies in particular that its fibers are contractible Kan
  complexes.
\end{proposition}

\begin{remark}
  Joyal proved the converse as well: if $\C$ is a simplicial set
  such that $\C^{\Delta^2} \to \C^{\Lambda^2_1}$ is a trivial Kan
  fibration, then $\C$ is a quasi-category.
\end{remark}

We can think of the map in the proposition roughly as, given a
``certified composition'', forgetting both the certificate and the
composite being certified. That the fibers of this map are
contractible says that to a homotopy theorist composition is uniquely
defined after all.

This result can be extended to strings of $n$ composable
$1$-simplices, namely, for a quasi-category $\C$ the canonical map
$\C^{\Delta^n} \to \C^{P_n}$ is a trivial Kan fibration. Here, $P_n =
\Delta^1 \vee_{\Delta^0} \Delta^1 \vee_{\Delta^0} \cdots
\vee_{\Delta^0} \Delta^1$ is the simplicial path of length $n$
obtained by glueing $n$ different $1$-simplices end to end. (When
$n=2$, it is isomorphic to $\Lambda^2_1$.) The case $n=3$ can be
interpreted as specifying a precise sense in which composition is
associative.

We can generalize even further to say that when defining a functor
from the free (ordinary) category on a directed graph $X$ into a
quasi-category $\C$, we can choose a diagram of $0$-simplices and
$1$-simplices in $\C$ of shape $X$ arbitrarily: there will always be
an extension to a functor, and moreover, the space of all such
extension is contractible. Formally, we have:

\begin{proposition}\label{freecats}
Let $X$ be a reflexive\footnote{Reflexive means
the graph has a distinguished loop at each vertex; these will play the
role of the identities in the free category on the graph.} directed
graph which we will think of as a simplicial\footnote{When thought of
a simplicial set, it is understood that the degenerate $1$-simplices
are the distinguished loops in the graph.} set which has no
non-degenerate $k$-simplices for $k \ge 2$. For any quasi-category
$\C$, the canonical map $\C^{NFX} \to C^X$ is a trivial Kan fibration,
where $NFX$ is the nerve of the free category on $X$.
\end{proposition}

For $X = P_n$, the free category on $X$ is the category which objects
$0,1,\ldots,n$ and a unique morphism from $i$ to $j$ when $i\le j$;
its nerve is the $n$-simplex $\Delta^n$, so we recover the previous
statement. As an example of this proposition, take $X$ to be a single
loop\footnote{Well, a single non-distinguished loop, in addition to
the distinguished one.}: an $X$-shaped diagram in a category is an
object together with an endomorphism. The free category on $X$ is just
the monoid of natural numbers under addition regarded as a one object
category, say, $FX = B\N$. In the world of ordinary categories, once
you've chosen an object and an endomorphism $f$ of it, you've uniquely
specified a functor out of $B\N$: the functor sends $k$ to $f^k = f
\circ f \circ \cdots \circ f$. For quasi-categories, there is no
canonical choice of $f^k$, you must make a choice for each $k$ and
then, to specify a functor out of $B\N$ you need to further choose
homotopies and higher homotopies showing you made compatible choices
of iterates of $f$. The proposition says then that all these choices
(of iterates and homotopies between their composites) can be made and
that, homotopically speaking, they are unique.

We haven't yet described how to get at mapping spaces in
quasi-categories. One intuitive approach is to use the arrow
\io-category of $\C$, which is simply the simplicial mapping space
$\C^{\Delta^1}$. This has a projection $\pi$ to $\C \times \C = \C^{\Delta^0
  \sqcup \Delta^0}$ which sends each $1$-simplex of $\C$ to its source
and target. Then, given two objects $X$ and $Y$ in $\C$, we can think
of them as being picked out by maps $\Delta^0 \to \C$ and form the
pullback:

\centerline{
  \xymatrix @R=5mm @C=3mm{
    \Map_\C(X,Y) \ar[d] \ar[r] & \C^{\Delta^1} \ar[d]^-{\pi} \\
    \Delta^0 \ar[r]_-{(X,Y)} & \C^{\Delta^0 \sqcup \Delta^0} \\
  }
}

This does work, that is, it produces a simplicial set $\Map_\C(X,Y)$
with the correct homotopy type, but there are many other descriptions
of the mapping spaces that are all homotopy equivalent but not
isomorphic as simplicial sets. One such alternative description of the
mapping spaces is given by Cordier's \emph{homotopy coherent nerve}
\cite{Cordier}, used in \cite{HTT} to compare quasi-categories with
simplicial categories. Cordier's construction not only provides models
for the mapping spaces but is also a procedure for strictifying
composition in a quasi-category: that is, constructing a simplicial
category (where composition \emph{is} required to be single-valued and
strictly associative) that represents the same \io-category as a given
quasi-category. Dan Dugger and David Spivak in \cite{DuggerSpivakRig}
explain a really nice way to visualize the mapping spaces appearing in
the homotopy coherent nerve through ``necklaces'' of simplices strung
together; they also wrote a second paper giving a detailed comparison
of the known constructions for mapping spaces in quasi-categories
\cite{DuggerSpivakMap}.

\subsection{Segal categories and complete Segal spaces}

Segal categories are a different formalization of the idea discussed
above for quasi-categories of a multivalued composition that is
uniquely defined homotopically. Just as quasi-categories can be
motivated by Proposition \ref{nerves}, Segal categories can be
motivated by the following equally easy result:

\begin{proposition}
  A simplicial set $X$ is isomorphic to the nerve of a category if and
  only if for each $n$, the canonical map $X_n \to X_1 \times_{X_0}
  X_1 \times_{X_0} \cdots \times_{X_0} X_1$ is a bijection.
\end{proposition}

This canonical map is the map $X^{\Delta^n} \to X^{P_n}$ we've already
met in section \ref{qcat}: it sends an $n$-simplex to its
\emph{spine}, the string of $1$-simplices connecting vertices $0$ and
$1$, $1$ and $2$, \dots, $n-1$ and $n$. It is tempting to try to make
compositions only defined up to homotopy simply by requiring these
canonical maps to be homotopy equivalences instead of bijections, but,
of course, that requires working with spaces rather than sets.

\begin{definition}
  A Segal category is a simplicial space (or more precisely a
  simplicial simplicial-set), that is, a functor
  $\DELTA^{\op} \to \sSet$ such that
  \begin{enumerate}
    \item the space of $0$-simplices $X_0$ is discrete, and
    \item for each $n$, the canonical map $X_n \to X_1 \times_{X_0}
          X_1 \times_{X_0} \cdots \times_{X_0} X_1$ is a weak homotopy
          equivalence.
  \end{enumerate}
\end{definition}

Complete Segal spaces, also called Rezk categories, were defined by
Charles Rezk in \cite{RezkCSS}; his purpose was explicitly
to find a nice model for the ``homotopy theory of homotopy theories'',
i.e., the \io-category of \io-categories\footnote{\io-categories
naturally form an $(\oo,2)$-category, but we can discard non-invertible
natural transformations to get an \io-category.}. Their definition is
a little complicated and we'll only describe it informally, but they
do have some advantages one of which was worked out by Clark Barwick
in his PhD thesis \cite{BarwickThesis}: the construction of complete
Segal spaces starting from simplicial sets as a model for
$\oo$-groupoids, can be iterated to provide a model for
$(\oo,n)$-categories. These are called \emph{$n$-fold complete Segal
spaces}, see \cite{ClarkChris} or \cite{LurieCob} for a definition, if
Barwick's thesis proves too hard to get a hold of.

A \emph{Segal space} like a Segal category, is also a simplicial
space, but we do not require that the space of objects $X_0$ be
discrete. In that case, the second condition must be modified to use
homotopy pullbacks\footnote{The reader unfamiliar with homotopy limits
  can find a quick introduction in section \ref{holim}.} so that it
reads: for each $n$, the canonical map $X_n \to X_1 \times^h_{X_0} X_1
\times^h_{X_0} \cdots \times^h_{X_0} X_1$ is a weak homotopy
equivalence. The \emph{completeness} condition has to do with the fact
that having a non-discrete space of objects means we have two
different notions of equivalence of objects: one is having an
invertible morphism between them in the \io-category modeled by $X$,
the other is being in the same connected component of $X_0$. Even
better, there are two canonical $\oo$-groupoids of objects: one is the
\emph{core} of the \io-category modeled by $X$, this is the
subcategory obtained by throwing away all non-invertible $1$-morphisms
(all higher morphisms are already invertible); the other is the
$\oo$-groupoid represented by $X_0$. The core of $X$ can be described
as a simplicial set directly in terms of the simplicial space $X$;
the completeness condition then says that it and $X_0$ are homotopy
equivalent.

\subsection{Relative categories}

Relative categories are based on the intuition that higher category
theory is meant for situations where we want to treat objects up to a
notion of equivalence that is weaker than isomorphism in the category
they live in. The reader should have in mind the examples of
equivalence of categories, Morita equivalence of rings, homotopy
equivalence of spaces, quasi-isomorphism of chain complexes, etc. The
definition of a relative category couldn't be simpler:

\begin{definition}
  A \emph{relative category} is a pair $(\C, \W)$ of an ordinary
  category $\C$ and a subcategory $\W$ of $\C$ required only to
  contain all the objects of $\C$. Morphisms in $\W$ are called
  \emph{weak equivalences}.
\end{definition}

Implicit in the claim that these somehow provide a model for
\io-categories is the claim that out of just a collection of weak
equivalences we get some sort of notion of homotopy between morphisms,
to play the role of $2$-morphisms in the \io-category represented by a
given relative category. To give the first idea of how this happens,
let's describe the homotopy category of the \io-category modeled by
$(\C,\W)$: it is $\C[\W^{-1}]$, the \emph{localization} of $\C$
obtained by formally adding inverses for all morphisms in $\W$. Let's
see in the example $\C = \Top$, $\W = \{$homotopy equivalences$\}$
that homotopic maps become equal\footnote{One's first instinct ---at
least, if one hasn't localized rings which are not integral domains---
might be that \emph{adding} inverses to some morphisms shouldn't force
other morphisms to become equal.} as morphisms in $\C[\W^{-1}]$.
First, notice that the projection $p : X \times [0,1] \to X$ is a
homotopy equivalence and thus becomes an isomorphism in
$\C[\W^{-1}]$. This means that the two maps $i_0, i_1 : X \to X\times
[0,1]$ given by $i_0(x) = (x,0)$ and $i_1(x) = (x,1)$ become equal in
the localization because $p \circ i_0 = p \circ i_1$. Finally, two maps
are homotopic when they can be written in the form $f \circ i_0$ and
$f \circ i_1$ for a single map $f$.

But of course, a satisfactory answer to the question of how higher
morphisms appear in the \io-category represented by $(\C, \W)$ would
construct the mapping space between two objects of $\C$, and this is
precisely what an enhancement of localization called \emph{simplicial
localization} does. We refer the reader to the classic papers of
William Dwyer and Dan Kan \cite{DK1}, \cite{DK2} for details on how
simplicial localizations may be constructed. The most intuitive
construction is probably the hammock localization, which we'll explain
by contrasting it with the (non-simplicial) localization
$\C[\W^{-1}]$.

Morphisms in $\C[\W^{-1}]$ can be represented by \emph{zig-zags}: $X_0
\xleftarrow{\sim} X_1 \to X_2 \xleftarrow{\sim} \cdots \to X_n$ where
arrows can go either way, but if they point to the left they are
required to be in $\W$ (this is typically indicated by placing a
$\sim$ on the arrow). \footnote{Really, the morphisms in $\C[\W^{-1}]$
  are equivalence classes of zig-zags in the smallest equivalence
  relation preserved by the operations of (1) removing an identity
  morphism, (2) composing two consecutive morphisms that point the
  same way, and (3) cancelling a pair of the form $\cdot
  \xrightarrow{w} \cdot \xleftarrow{w} \cdot$ or $\cdot \xleftarrow{w}
  \cdot \xrightarrow{w} \cdot$.} To form the mapping space in the
hammock localization we add homotopies between zig-zags: the mapping
space is constructed as a simplicial sets whose vertices are zig-zags
and whose $1$-simplices are weak equivalences of zig-zags, by which we
mean diagrams of the form:

\centerline{ \xymatrix @R=1.5mm @C=5mm { & X_1\ar[dd]^[@]\sim & X_2
\ar[dd]^[@]\sim\ar[l]_\sim\ar[r] & X_3 \ar[dd]^[@]\sim &
\cdots\ar[l]_\sim\ar[r] & X_{n-1} \ar[dd]^[@]\sim \ar[dr] & \\ X_0
\ar[ur]\ar[dr] & & & & & & X_n \\ & Y_1 & Y_2 \ar[l]_\sim\ar[r] & Y_3
& \cdots\ar[l]_\sim\ar[r] & Y_{n-1} \ar[ur] & \\ } }

\noindent in which all the left pointing morphisms and all the
vertical ones are required to be in $\W$. Higher dimensional simplices
are similar but have more rows (and look even more like hammocks than
$1$-simplices do).\footnote{In this description some details are
missing, see \cite{DK2}.}

\begin{further}
The papers by Dwyer and Kan on simplicial localization
already indicate that relative categories, bare-bones though they may
be, can be used to model \io-categories. The book \cite{DKHS} develops
homotopy theory for relative categories (there called ``homotopical
categories'' and required to satisfy a mild further axiom).
More recently, Clark Barwick and Dan Kan, in a series of papers
\cite{BK1}, \cite{BK2}, \cite{BK3}, compare relative categories
to other models of \io-categories and
define a generalization of them that provides a model for
$(\oo,n)$-categories.
\end{further}

\subsection{$A_\oo$-categories}\label{top-Aoo}

An \emph{operad}\footnote{Technically, a non-symmetric operad
as we don't ask for an action of the symmetric groups.} is a
collection of spaces $\O(n)$ together with composition maps
\[ \O(n) \times \O(k_1) \times \O(k_2) \times \cdots \times
     \O(k_n) \to \O(k_1 + k_2 + \cdots + k_n), \]
which are required to satisfy associativity conditions\footnote{We
are also omitting a couple of conditions on $\O(0)$ and $\O(1)$.}
that are easy to guess if one thinks of the process of substituting
$n$ functions, one of $k_1$ variables, one of $k_2$ variables, etc.,
into a function of $n$ variables to obtain an overall function of
$k_1 + k_2 + \cdots + k_n$ variables. As this suggests, the
elements of $\O(n)$ are called $n$-ary operations. They can
be used to parametrize all the homotopies required for a composition
that is associative up to coherent homotopy:

\begin{definition}
  An $A_\infty$-operad is one such that all $\O(n)$ are contractible.
  Given any such operad, an $A_\infty$-category $\C$ consists of
  \begin{enumerate}
    \item a collection of objects,
    \item a space $\C(X,Y)$ for every pair of objects, and
    \item composition maps
        \[\O(n) \times \C(X_{n-1}, X_n) \times \C(X_{n-2}, X_{n-1})
        \times \cdots \times \C(X_0, X_1) \to \C(X_0, X_n)\] for every $n$
        and every sequence of objects $X_0$, $X_1$, \dots, $X_n$,
  \end{enumerate}
  which are required to be compatible with the composition operations
  of the operad in an obvious sense.\footnote{In both the description
  of operad and $A_\infty$-category we've
  omitted discussing identities. The reader can easily supply the
  missing details.}
\end{definition}

Notice, just like topological categories, this model provides easy access
to homotopy categories: since $\O(n)$ is contractible, applying $\pi_0$
to an $A_\oo$-category produces an ordinary category with $\Hom$-sets
given by $\pi_0(\C(X,Y))$.

In case the $A_\oo$-operad $\O$ is the operad of Stasheff associahedra,
an $A_\oo$-category with a single object is equivalent to the original
notion of an $A_\oo$-space introduced by Jim Stasheff in his work
on homotopy associative $H$-spaces \cite{Stasheff}.
The main result of that work can be interpreted as proving the
homotopy hypothesis for $A_\oo$-categories with a single object. We'll
state a less precise version informally:

\begin{proposition}
  An $A_\infty$-space $X$ is weak homotopy equivalent to a loop space
  $\Omega Y$ (in such a way that composition in $X$ corresponds to
  concatenation of loops) if and only if $\pi_0(X)$ is a group under
  the operation induced from composition in $X$.
\end{proposition}

This \emph{recognition principle} for loop spaces is part of the
original motivation for Peter May's definition of operad\footnote{
  Which, we repeat, besides the data mentioned above includes actions
  of $\Sigma_n$ on $\O(n)$ and requires the composition maps to be
  equivariant in an easily guessed sense.}, which he used to prove
a recognition principle in the same spirit for iterated loop spaces
$\Omega^n Y$ and infinite loop spaces (which can be thought of as a
sequence of spaces $Y_0, Y_1, \ldots$ each of which is equivalent to
the loop spaces of the following one). See May's book
\cite{MayGILS}.

\begin{further}
  This model doesn't seem to get used that much in practice. The only
  example of a paper constructing some \io-category as an
  $A_\oo$-category that the author is aware of is \cite[Proposition
  1.4]{Ching}. Todd Trimble used $A_\oo$-categories as the first step
  in an inductive definition of $(\oo,n)$-category, see
  \cite{ChengOperadic}. A talk given in Morelia by Peter May
  \cite{MayMorelia} expressed the hope that a simpler, more general
  version of the inductive approach would work, but it was pointed out
  by Michael Batanin that this doesn't quite work (this is mentioned
  in Eugenia Cheng's paper just cited).
\end{further}

\subsection{Models of subclasses of \io-categories}

There are also several ways of modeling special classes of
\io-categories, which, when applicable can be simpler to calculate
with. We'll mention model categories and derivators for which it is
hard to say exactly which \io-categories they can model, but which
certainly can only model \io-categories that have all small homotopy
limits and colimits, and linear models, which model \io-categories
that are enriched over an \io-categories of chain complexes.

\subsubsection{Model categories}

Quillen's model categories are the most successful setting for
abstract homotopy theory. A \emph{model category} $\C$ is an ordinary
category that has all small limits and colimits and is equipped with
three collections of morphisms called \emph{cofibrations},
\emph{fibrations} and \emph{weak equivalences} which are required to
satisfy axioms that abstract properties that hold of the classes of
maps of topological spaces that they are named after. We won't give a
precise definition, but refer the reader to standard references such
as the introduction \cite{DwyerSpalinski}, the book \cite{Hovey}, the
fast paced \cite[Appendix A.2]{HTT} or Quillen's original
\cite{Quillen} (but note that what we are calling model categories
where called \emph{closed} model categories there). The \io-category
modeled by a model category $\C$ is just the one modeled by the
relative category $(\C, \text{weak equivalences})$. We mention them
separately from relative categories because the extra structure makes
them easier to deal with than a random relative category, so they form
an eminently practical way to construct particular \io-categories.
Even on the level of homotopy categories, the fibrations, cofibrations
and the axioms make the localization better behaved. For example, in
the homotopy category of a model category we do not need to consider
zig-zags of arbitrary length, it is enough to look at zig-zags of the
form $\cdot \xleftarrow{\sim} \cdot \to \cdot \xleftarrow{\sim}$.

\begin{remark}
  Small homotopy limits and colimits always exist in a model category
  and thus they can only model \io-categories that are complete and
  cocomplete. It is not known to the author whether or not all such
  \io-categories arise from model categories. There is however a
  result of Carlos Simpson's under further smallness assumptions,
  namely he showed that \emph{combinatorial} model categories provide
  models precisely for the class of \emph{locally presentable}
  \io-categories.  See \cite{SimpsonPres} (but beware that what are
  called cofibrantly generated model categories there are what we are
  calling combinatorial model categories), or \cite[Section
  5.5.1]{HTT}. Roughly speaking, a locally presentable category is one
  that is cocomplete and generated under colimits by a small
  subcategory of objects which are small or compact in some sense. A
  combinatorial model category is required to be locally presentable
  (and to have a model structure which is cofibrantly generated, which
  is also a condition with the flavor of the whole being determined by
  a small portion). For information about Jeff Smith's notion of
  combinatorial model category see \cite[Appendix A.2.6]{HTT}.
\end{remark}

Model categories have been hugely successful in providing workable
notions of homotopy theory in many topological and algebraic contexts.
A wealth of model structures have been constructed and all provide
examples of \io-categories that people care about. When performing
further constructions based on these \io-categories, such as taking
categories of diagrams in one of them, functors between two of them or
homotopy limits or colimits of them it can be very hard to remain in
the world of model categories. In those cases, using model categories
to present the inputs to these constructions but carrying them out in
the world of \io-categories is a very reasonable compromise.

We will meet model categories again in section \ref{comparison} since
the original comparison results between models of \io-categories were
formulated in that language.

\subsubsection{Derivators}

When working with an \io-category $\C$, it might be tempting to do as
much as possible in $\ho\C$, since ordinary categories are much
simpler and more familiar objects. We can't get very far, though, we
run into trouble as soon as we start talking about homotopy limits and
colimits\footnote{See section \ref{holim} for a quick introduction.}.
Say we have a small (ordinary) category $\I$ and wish to talk about
homotopy limits or colimits of $\I$-shaped diagrams in $\C$. Homotopy
limits should be homotopy invariant: if two diagrams $F,G : \I \to \C$
are connected by a natural isomorphism\footnote{A natural
transformation whose components are invertible in the sense we always
use for \io-categories: invertible up to higher invertible
morphisms.}, they should have equivalent limits in $\C$. So, taking
homotopy limits should induce a functor $\ho(\C^\I) \to \ho\C$. Now,
this homotopy category $\ho(C^\I)$ is not something we can construct
just from $\ho \C$ and $\I$, in particular it is not equivalent to
$(\ho\C)^\I$.

\begin{example} Let $\I$ be $\Z/2$ regarded as a category with a
single object and let $\C$ be the \io-category of $\oo$-groupoids (or
spaces). An $\I$-shaped diagram in $\C$ is just a space with an action
of $\Z/2$. Consider the diagrams given by the trivial action and the
$180^\circ$ rotation on $S^1$. Since the $180^\circ$ rotation is
homotopic to the identity on $S^1$, these two diagrams become equal in
$(\ho\C)^\I$, but are not isomorphic in $\ho(\C^\I)$ since, for
example, they have different homotopy colimits: since the rotation
action is free, the homotopy colimit in that case is just $S^1/(\Z/2)
\cong S^1$; for the trivial action, we get $(E(\Z/2) \times
S^1)/(\text{diagonal action}) = B(\Z/2) \times S^1$. \end{example}

The idea of derivators then, is to hold on to, not just $\ho\C$, but
$\ho(\C^\I)$ for every small (ordinary) category $\I$ as well. This at
least allows one to hope to be able to discuss homotopy limits and
colimits. Given an \io-category $\C$, the construction $\I \to
\ho(\C^\I)$ provides a strict $2$-functor $(\Cat_\text{small})^{\op}
\to \Cat$ where $\Cat$ is the strict $2$-category of all not
necessarily small categories and $\Cat_\text{small}$ is the
sub-$2$-category of small ones.\footnote{If the reader is not versed
in the art of worrying about size issues, we advise not to start until
after reading this survey. We do however caution that while it might
seem like a merely technical point there is substance to it: for
example, it is easy to prove that if the collection of morphisms of a
category has size $\lambda$ and the category has products
$\prod_{i<\lambda} X_i$ then it is a preorder.} By definition,
derivators are strict $2$-functors $(\Cat_\text{small})^{\op} \to
\Cat$ satisfying further conditions that guarantee that homotopy
limits and colimits (and more generally homotopy versions of the left
and right Kan extensions) exist and are well-behaved. As in the case
of model categories: (1) the definition directly implies derivators
can only model \io-categories which are complete and cocomplete, (2)
the author does not know if all such \io-categories can be modeled,
and (3) adding presentability on both sides of the equation balances
it, see \cite{RenaudinPres}. The later \cite{RenaudinLP} deals with
representing \io-categories coming from left proper model categories
by derivators.

\begin{further} Derivators were defined by Alexander Grothendieck (the
term appears first in \cite{PursuingStacks}, a few years later
Grothendieck wrote \cite{Derivateurs}) and independently by Alex
Heller \cite{Heller} (who called them ``homotopy theories''). Good
introductions can be found in \cite{Maltsiniotis}, \cite{GrothDer}, and
the review section of \cite{ShulmanPontoGroth}. \end{further}

\subsubsection{dg-categories, $A_\oo$-categories}

Now we'll discuss two ``linear'' models (or more precisely, models
based on chain complexes) for special kinds of \io-categories that
have seen much use in algebra and algebraic geometry. These are the
notions of dg-categories, which are analogous to topological or
simplicial categories, and $A_\oo$-categories, which are analogous to
the identically named $A_\oo$-categories mentioned in the previous
section (these chain complex based $A_\oo$-categories see much more
use than their topological counterparts and most people associate the
name $A_\oo$-category with the chain complex version described in this
section). In both cases the analogy comes from replacing spaces by
chain complexes, that is, by restricting $\oo$-groupoids to those
modeled by chain complexes: the abelian, fully strict $\oo$-groupoids.
(Recall the Dold-Kan correspondence which establishes an equivalence
of categories between simplicial abelian groups and chain complexes of
abelian groups concentrated in non-negative degrees.)

\begin{definition} Let $R$ be a commutative ring. A \emph{differential
graded category} or \emph{dg-category} over $R$ is a category enriched
in the monoidal category of chain complexes of $R$-modules.
\end{definition}

(There are variants where the complexes are required to be bounded,
and the analogy discussed above makes the most sense directly for
complexes of abelian groups in non-negative degrees.)

$A_\oo$-categories can be defined using an operad in chain complexes,
analogously to the definition in section \ref{top-Aoo}, but there is a
more explicit description of them in terms of one $n$-ary composition
operation for each $n$. We will not repeat the definition here, since
we won't have anything to say about it, but the interested reader can
look it up in the references.

\begin{further} Good introductions to dg-categories include
\cite{Keller} and \cite{ToenDG}. The \io-category of all dg-categories
is described as a model category in Gon\c{c}alo Tabuada's PhD thesis
\cite{Tabuada}. That model structure is used in \cite{ToenMorita}
to develop Morita theory for dg-categories.
A concise introduction to $A_\oo$-categories can be
found in \cite{Keller}; a thorough reference is the book
\cite{Bespalov}. Also see Maxim Kontsevich and Yan Soibelman's notes
\cite{KontsevichSoibelman} which deal mostly with $A_\oo$-algebras
which are $A_\oo$-categories with a single object (although that
definition is anachronistic, of course). The most conspicuous example
of an $A_\oo$-category is the Fukaya category of a symplectic manifold
that plays a starring role in homological mirror symmetry. It was
first mentioned in \cite{Fukaya}. For comprehensive treatments see the
books \cite{FOOO} and \cite{Seidel}, the later of which has a good
introduction to $A_\oo$-categories in Chapter 1. \end{further}

\section{The comparison problem}\label{comparison}

Comparing different definitions of higher categories is harder than it
might seem at first, since even what it means to show that two
theories are equivalent is not completely clear. Imagine we wish to
compare two theories of $n$-categories, call
them red categories and blue categories. Just like ordinary categories
form a $2$-category, red $n$-categories should form an
$(n+1)$-category, and we'd like to show that this $(n+1)$-category
is equivalent in an appropriate sense to the $(n+1)$-category of
blue $n$-categories. But what color are the $(n+1)$-categories and the
equivalence? It is reasonable to expect that red $n$-categories form
a red $(n+1)$-category and similarly that the blue $n$-categories
form a blue $(n+1)$-category. But this means we can't easily compare
the $(n+1)$-categories before solving the comparison for red and blue
higher categories!

We might be able to assemble all blue $n$-categories into a red
$(n+1)$-category by an ad hoc construction, and show that that red
$(n+1)$-category is red-equivalent to the red $(n+1)$-category of
red $n$-categories. In that case we could say the red theory regards
the two theories as equivalent. But a priori, if that happens, the blue
theory might disagree and not consider the two theories equivalent!

The way the comparison problem was solved for \io-categories is as
follows:

\begin{enumerate}
  \item What gets compared are not the $(\oo,2)$-categories of all
    \io-categories of particular kind, but rather \io-categories of
    \io-categories (obtained from the $(\oo,2)$-category by throwing
    away non-invertible natural transformations).
  \item For each model $M$ the \io-category $\Catio^M$ of all
    $M$-style \io-categories was described as a model category.
    \footnote{This does not necessarily mean that a model structure
      was put precisely on some ordinary category category of all
      $M$-style \io-categories; but rather on a larger category in
      which the objects which are both fibrant and cofibrant are the
      $M$-style \io-categories. For example, for quasi-categories the
      model structure is on the category of all simplicial sets: every
      object is cofibrant and the fibrant ones are precisely the
      quasi-categories.}
  \item It was shown that these model categories are all connected by
    zig-zags of \emph{Quillen equivalences}: these are equivalences
    of categories that preserve enough aspects of the model structure
    to ensure that two model categories have the same homotopy theory,
    i.e., model the same \io-category.
\end{enumerate}

In the non-standard terminology above, this says that the theory of model
categories regards all models as equivalent. It can also be shown that
any model regards all models as equivalent. For example, take quasi-categories.
One can construct a quasi-category from a model category and show that
two model categories connected by a zig-zag of Quillen equivalences produce
quasi-categories that are equivalent (according to a particular definition
of equivalence of quasi-categories). Then we get a quasi-category $\Catio^M$
for each of the models, and they are all equivalent.

Portions of the program outlined above showing the equivalence of five
of the models we discussed (namely quasi-categories, simplicial
categories, Segal categories, complete Segal spaces and relative
categories) were carried out by Julie Bergner \cite{BergnerEquiv},
Clark Barwick and Dan Kan \cite{BK2}, \cite{BK2}, Andr\'{e} Joyal and
Miles Tierney \cite{JoyalNotes}, \cite{JoyalTierney}. For a beautiful
diagram showing the Quillen equivalences at a glance and further
references see \cite[Figure 1]{ClarkChris}. For a summary of the model
structures and Quillen equivalences comparing the first four models
(i.e., excluding relative categories) see Julie Bergner's survey
\cite{BergnerSurvey}.

\subsection{Axiomatization}

A recent breakthrough in the theory of $(\oo,n)$-cat\-e\-go\-ries is
the axiomatization by Clark Barwick and Chris Schommer-Pries
\cite{ClarkChris} of the \io-category $\Catin$ of
$(\oo,n)$-categories. As in the direct comparison results we mentioned
for \io-categories, what gets axiomatized is not an
$(\oo,n+1)$-category  of $(\oo,n)$-categories but rather the
\io-category one obtains from that $(\oo,n+1)$-category by throwing
away higher non-invertible morphisms.  Their work was inspired by
Bertrand To\"{e}n's influential \cite{ToenAx} that similarly
axiomatizes \emph{model categories} of \io-categories.  To\"{e}n's
axioms are closely related to Giraud's axioms for toposes, while
Barwick and Schommer-Pries's axioms stray a bit further from them.
Very roughly, Barwick and Schommer-Pries axioms are as follows:

\begin{enumerate}
  \item There is an embedding in $\Catin$ of the category of
    \emph{gaunt} $n$-categories. A gaunt $n$-category is a strict
    $n$-category all of whose invertible $k$-morphisms are identities,
    for all $k$. The images of these gaunt $n$-categories are
    required to generate $\Catin$ under homotopy colimits.
  \item The embedding provides us in particular with
    $(\oo,n)$-categories that are ``walking $k$-cells'': they consist
    of a single $k$-morphism and its required (globular) boundary.
    Certain gaunt $n$-categories obtained by glueing a few cells
    together are required to still be obtained by the same glueing
    process inside $\Catin$, that is, the embedding of gaunt
    $n$-categories is required to preserve a few colimits of diagrams
    of cells.
  \item $\Catin$ is required to have internal $\Hom$s, and more
    generally, slices $\Catin/C_k$ are required to have internal
    $\Hom$s for all cells $C_k$.
  \item The \io-category $\Catin$ is \emph{universal} with the respect
    to the above properties, in the sense that any other \io-category
    $\D$ with an embedding of gaunt $n$-categories satisfying the first
    three axioms is obtained from $\Catin$ via a localization functor
    $\Catin \to \D$ that commutes with the embeddings.
\end{enumerate}

Even more roughly: the first axiom gives us the ability to present
$(\oo,n)$-categories by means of generators and relations, and the
second guarantees that at least for presentations of gaunt
$n$-categories we get the expected answer. There are many imprecisions
in our description of the axioms and we refer the reader to
\cite{ClarkChris} for a correct statement.

\begin{remark}
  One might think that discarding those morphisms loses too much
  information, but in a sense it doesn't: the \io-category
  $\Catin$ characterized by the axioms is Cartesian closed, so
  for any two $(\oo,n)$-categories $\C$ and $\D$ there is an
  $(\oo,n)$-category $\Fun(\C,\D)$ that contains all the non-invertible
  natural transformations and higher morphisms that are not directly
  observable in the mapping space $\Map_{\Catin}(\C,\D)$.
\end{remark}

To\"{e}n also proved in \cite{ToenAx} that the only automorphism of
the theory of \io-categories is given by taking opposite categories.
It is unique in the strong sense that not only is every automorphism
equivalent to that one, but automorphisms themselves have no
automorphisms, or more precisely, there is a naturally defined
$\oo$-groupoid of automorphisms of $\Catio$ and it is homotopy
equivalent to the discrete group $\Z/2$. Barwick and Schommer-Pries
prove the analogue of this result for $\Catin$, showing that its
space of automorphisms is the discrete group $(\Z/2)^n$ corresponding
to the possibility of deciding for each degree $k$ separately whether
or not to flip the direction of the $k$-morphisms.

\section{Basic \io-category theory}

In this section we review how the most basic concepts of category
theory generalize to \io-categories. The philosophy is that
\io-categories are much more like ordinary 1-categories than fully
weak $\omega$-categories where no morphisms are required to be
invertible\footnote{In fact, even $2$-categories have some tricky
  points that that do not arise when dealing with \io-categories, such
  as the need to distinguish between several kinds of limits and
  colimits called pseudo-limits, lax-limits and colax-limits.}. A
little more precisely, the intuition that \io-categories have
spaces of morphisms and that these spaces only matter up to (weak)
homotopy equivalence usually leads to useful definitions and
correct statements. We will also frequently point out what happens in
the case of quasi-categories, which due both intrinsically to some
features they possess and externally to the availability of
\cite{Joyal}, \cite{HTT} and \cite{HA} are probably the most
``practical'' model to work with.

\begin{further} Besides those systematic treatises by Joyal and Lurie
already mentioned, we recommend \cite{Groth}, an excellent and well
motivated summary of large chunks of \cite{HTT}, \cite{DAG1},
\cite{DAG2} and \cite{DAG3} (the last three of which were reworked
into the first few chapters of \cite{HA}). The availability of
\cite{Groth} is why we feel justified in giving very little detail in
this section, just giving the flavor of the topic. Also highly
recommended is the forthcoming book \cite{Riehl}, Part IV of which is
about quasi-categories and includes, among other things, (1) a
discussion of which aspects of \io-categories are already captured by
the $2$-category whose objects are quasi-categories, whose morphisms
are functors and whose $2$-morphisms are homotopy classes of natural
transformations; and (2) plenty of geometrical information about
quasi-categories viewed as simplicial sets, such as how to visualize
mapping spaces in quasi-categories or the homotopy coherent nerve that
produces the equivalence between quasi-categories and simplicial
categories. \end{further}

\subsection{Equivalences}

There are two things typically called \emph{equivalences} in
\io-category theory: one generalizes isomorphisms in an ordinary
category, and the other generalizes equivalences between categories.
These are related in the following way: we can ignore
non-invertible natural transformations to get $\Catio$, the
\io-category of \io-categories, functors and invertible natural
transformations. A functor is then an equivalence of \io-categories
if and only if it is an equivalence as a morphism in $\Catio$.

\begin{remark}
  That \io-isomorphisms go by the name equivalences is probably due
  to (1) the case of the \io-category of spaces where they are just
  (weak) homotopy equivalences, and (2) there being, in the model of
  topological categories (resp. simplicial categories), a second, stricter
  notion of isomorphism, coming from enriched category theory:
  morphisms such that composition with them induces homeomorphisms
  (resp. isomorphisms of simplicial sets) between mapping spaces.
\end{remark}

\begin{definition}
  A morphism $f : X \to Y$ in an \io-category is an \emph{equivalence}
  if its image in $\ho\C$ is an isomorphism, or, equivalently, if for
  every object $Z \in \C$, $f \circ \_ : \Map_\C(Z,X) \to \Map_\C(Z,Y)$
  and $\_ \circ f : \Map_\C(Y,Z) \to \Map_\C(X,Z)$ are weak homotopy
  equivalences.
\end{definition}

Andr\'{e} Joyal proved, using the model of quasi-categories, that when
a morphism $f$ is an equivalence, one can coherently choose an inverse
$g$ for it, $2$-morphisms showing $f$ and $g$ are inverses, $3$-morphisms
showing invertibility of those $2$-morphisms, and so on. More precisely,
such a coherent system of choices is given by a functor $F : \J \to \C$
where $\J$ is the ordinary category $0 \xrightarrow{\simeq} 1$,
which has two objects and a unique isomorphism between them. The precise
statement then is as follows:

\begin{proposition}[{\cite[Corollary 1.6]{JoyalPub}}]
  A morphism $f$ in an \io-category $\C$ is an equivalence if and only
  if there is a functor $F : \J \to \C$ with $F(0 \to 1) = f$, where
  $\J$ is the \emph{walking isomorphism} defined above.
\end{proposition}

This result is in the same spirit as the much easier automatic
homotopy coherence result in \ref{freecats}.

\subsubsection{Further results for quasi-categories} The idea behind
requiring the horn filler condition only for inner horns, $\Lambda^n_k$
with $0<k<n$ in the definition of quasi-category is that, for $n=2$,
filling a $\Lambda^2_0$ horn requires inverting the edge between vertices
$0$ and $1$:

\centerline{$ \begin{gathered} \xymatrix @R=6mm @C=2mm {
& Y & \\ X \ar[ur]^{f} \ar[rr]_{g} & & Z } \end{gathered}
\hspace{7mm}\leadsto\hspace{7mm} \begin{gathered} \xymatrix @R=6mm
@C=2mm { & Y \ar[dr]^{h \simeq g \circ f^{-1}} \ar@{=>}[d]_-{\alpha}& \\ X \ar[ur]^{f}
\ar[rr]_{g} & & Z } \end{gathered}$ }

And indeed Joyal proved that this is the only obstacle:

\begin{proposition}[{\cite[Theorem 1.3]{JoyalPub}}]
  In a quasi-category $\C$, a morphism $f$ is an equivalence if and
  only if there are fillers for every horn $\Lambda^n_0 \to \C$
  whose edge joining vertices $0$ and $1$ is $f$.
\end{proposition}

There is, of course, a dual result about $\Lambda^n_n$ horns. Putting
the two together we get the homotopy hypothesis for quasi-categories:
every $1$-simplex of a quasi-category $\C$ is an equivalence if and
only if $\C$ is a Kan complex.

\subsection{Limits and colimits}\label{holim}

The notion of limit and colimit in an \io-category should be thought
of as generalizations of \emph{homotopy} limits and colimits, and
indeed, reduce to those for the \io-category of spaces (or more
generally, for an \io-category coming from a model category; there
is also a notion of homotopy limit and colimit in a model category).
We begin by recalling what homotopy limits and colimits of
diagrams of spaces are, by describing a construction of them that
provides the correct intuition for \io-categories. When reading
the following description, the reader should think of paths as
\emph{invertible} morphisms in an \io-categories, and homotopies
as higher morphisms.

Let's consider a diagram of spaces $F : \I \to \Top$ (where, for now,
$\I$ is an ordinary category). Recall that its (ordinary,
non-homotopy) limit can be constructed as follows: it consists of the
subspace of $\prod_{i \in \I} F(i)$ of points $(x_i)_{i\in\I}$ such
that for any $\alpha : i \to j$ in $\I$, we have that $F(\alpha)(x_i)
= x_j$. To construct instead the homotopy limit, following the
philosophy described in remark \ref{metaphor}, we replace the strict
notion of equality in $F(\alpha)(x_i) = x_j$ by the corresponding
homotopical notion: a path from $F(\alpha)(x_i)$ to $x_j$, this path
is witness to the fact that ``homotopically speaking, $F(\alpha)(x_i)$
and $x_j$ are the same''. For simple diagrams, without composable
arrows, this is enough.  For example, the homotopy pullback of a
diagram $X \xrightarrow{f} Z \xleftarrow{g} Y$ \emph{can} be
constructed as the subspace of $X \times Y \times Z \times Z^{[0,1]}
\times Z^{[0,1]}$ consisting of $5$-tuples $(x,y,z,\gamma,\sigma)$
with $\gamma(0) = f(x), \gamma(1) = z = \sigma(1), \sigma(0) = g(y)$.

For diagrams that have pairs of composable arrows, we need to ensure
these paths that act as witnesses compose as well, that is, if we have
$i \xrightarrow{\alpha} j \xrightarrow{\beta} k$ somewhere in $\I$,
every point of the homotopy limit will have among its coordinates the
data of
\begin{enumerate}
  \item points $x_i \in F(i), x_j \in F(j), x_k \in F(k)$, and
  \item paths $\gamma_{ij} : F(\alpha)(x_i) \to x_j$,
   $\gamma_{jk} : F(\beta)(x_j) \to x_k$ and $\gamma_{ik} :
   F(\beta\circ\alpha)(x_i) \to x_j$.
\end{enumerate}
The paths $\gamma_{ik}$ and $\gamma_{jk} \cdot
(F(\beta)\circ\gamma_{ij})$ both witness that $x_i$ and $x_j$ are in
the same path component of $F(k)$, but we shouldn't regard their
testimonies as being independent! The points of the homotopy limit
should also include a homotopy between these two paths. Clearly, this
doesn't stop here: for diagrams with triples of composable arrows we
should have homotopies between homotopies and so on.  For any given
small diagram shape $\I$ it is clear which paths and homotopies are
required.

Dually, for homotopy colimits, instead of taking a subspace of a large
product where we require the presence of some paths and homotopies, we
form a large coproduct of spaces and glue in paths and homotopies
that enforce ``sameness'' of points in the colimits.

\begin{remark}
  There is also a general formula for the required homotopies in
  homotopy limits due to Pete Bousfield and Dan Kan
  \cite{BousfieldKan}: for a functor $F : \I \to \Top$, we have
  \[\holim F = \int_{i\in\I} F(i)^{\abs{N(I\downarrow i)}}.\] Here
  $\abs{N(I\downarrow i)}$ is the geometric realization of the nerve
  of the slice category $\I \downarrow i$ whose objects are objects of
  $\I$ with a map to $i$, and whose morphisms are commuting triangles.
  The $\int$ sign denotes a type of limit called an \emph{end}, and
  can be constructed as a subspace of the product of function spaces
  $\prod_{i \in \I} F(i)^{\abs{N(\I\downarrow i)}}$ consisting of
  compatible families of functions $\gamma_i :\abs{N(\I\downarrow i)}
  \to F(i)$: a morphism $i \xrightarrow{\alpha} j$ induces a functor
  $\I\downarrow i \to \I\downarrow j$ by composition, and thus gives a
  map $\alpha_* : \abs{N(\I\downarrow i)} \to \abs{N(\I\downarrow
    j)}$; the family $\{\gamma_i\}$ is compatible if $\gamma_j \circ
  \alpha_* = F(\alpha) \circ \gamma_i$. The reader can check that the
  homotopies mentioned above for composable pairs
  $i\xrightarrow{\alpha} j\xrightarrow{\beta} k$ are encoded by this
  formula as maps out of triangles that restrict on the boundary to
  the paths corresponding to the morphisms $\alpha$, $\beta$ and
  $\beta \circ \alpha$.
\end{remark}

All limits and colimits in \io-categories can be defined in terms of
homotopy limits of spaces. Recall that for a diagram $F : \I \to \C$
in an ordinary category $\C$, the limit and colimit can be defined by
requiring the canonical functions of sets
\[\Hom_\C(X,\lim_{i\in\I} F(i)) \to \lim_{i\in\I} \Hom_\C(X,F(i))\] and
\[\Hom_\C(\colim_{i\in\I} F(i),X) \to \lim_{i\in\I^{\op}}
  \Hom_\C(F(i),X)\]
to be bijections natural in $X$. The limits occurring in the codomain
of these canonical maps are taken in the category of
sets\footnote{Note that \emph{co}limits in the category of sets do not
  appear here.}.

\begin{definition}
  Given a functor $F : \I \to \C$ between two \io-categories, its
  \emph{limit} and \emph{colimit}, if they exist, are determined
  up to equivalence in $\C$ by requiring that
  \[\Map_\C(X,\lim_{i\in\I} F(i)) \to \holim_{i\in\I} \Map_\C(X,F(i))\] and
  \[\Map_\C(\colim_{i\in\I} F(i),X) \to \holim_{i\in\I^{\op}} \Map_\C(F(i),X)\]
  are weak equivalences of spaces, natural in $X$.
\end{definition}

Note the special case of initial and terminal objects: an object $X$ of
an \io-category $\C$ is initial if $\Map_\C(X,Y)$ is contractible for all
objects $Y$, and terminal if $\Map_\C(Y,X)$ is contractible for all
$Y$. As expected initial objects are unique when they exist: more
precisely, Joyal proved that the full subcategory of $\C$ consisting
of all initial objects is a contractible $\oo$-groupoid, that is, any
two initial objects are equivalent, any two equivalences between
initial objects are homotopic, etc.

Other constructions of limits and colimits in ordinary categories also
generalize to the \io-setting and give the same definition as above.
For example, the limit of a functor $F$ can be characterized as a
terminal object in the category of cones over $F$. For concreteness
we'll use quasi-categories to describe how this goes for
\io-categories: given a functor $F : \I \to \C$, the quasi-category of
cones over $F$, $\Cones(F)$, is the simplicial set whose $n$-simplices
are given by maps of simplicial sets $\Delta^n \star \I \to \C$ which
restrict to $F$ on $\I$. Its vertices are exactly what we'd expect,
since $\Delta^0 \star \I$ is $\I$ with a new initial object adjoined.
Here $K \star L$ denotes the join of simplicial sets, a geometric
operation that can be thought as providing a canonical triangulation
of the union of line segments joining all pairs of a point of $\abs{K}$
and a point of $\abs{L}$; each $k$-simplex of $\I$ contributes an
$(n+k+1)$-simplex to $\Delta^n \star \I$.

\begin{remark}
  While it is perfectly fine to define a diagram in a category to be a
  functor for developing the theory, most people don't actually think
  of many common diagram shapes, such as diagrams for pullbacks,
  pushouts, infinite sequences $X_0 \to X_1 \to \cdots$, etc., as
  being given by a category; instead these shapes are usually thought
  of as being given by a directed graph. Another advantage of the
  model of quasi-categories is that there is a very convenient
  generalization of directed graphs and diagrams shaped like them:
  simply take arbitrary simplicial sets for shapes and define a
  $K$-shaped diagram in $\C$ to be a map of simplicial sets $K \to
  \C$. Both descriptions given above make sense for these more general
  types of diagrams. Also, this notion of $K$-shaped diagram can be
  used in other model of \io-categories but is more cumbersome, since
  essentially one needs to define the free \io-category on $K$.
\end{remark}

Most classical results about limits and colimits hold for
\io-categories with appropriate definitions, and we hope that the
examples shown here give a rough idea of how the definitions are
generalized.

\subsection{Adjunctions, monads and comonads}\label{adj}

As for other concepts, the definition of adjunction in ordinary
category theory that uses $\Hom$-sets and bijections generalizes to
\io-categories by using mapping spaces and homotopy equivalences:

\begin{definition}
  Given functors $F : \C \to \D$ and $G : \D \to \C$ between
  \io-categories, an \emph{adjunction} is specified by a giving a
  \emph{unit}, a natural transformation\footnote{We haven't defined
    natural transformations between functors of \io-categories. They
    are morphisms in functor \io-categories; here, $u$ is a morphism
    in $\Fun(\C,\C)$.} $u : \id_\C \to g \circ f$
  such that the composite map
  \[ \Map_\D(F(C),D) \xrightarrow{G} \Map_\C(G(F(C)),G(D))
    \xrightarrow{\_ \circ u_C} \Map_\C(C,G(D)) \]
  is a weak homotopy equivalence.
\end{definition}

As in the case of ordinary categories, if $F$ has a right adjoint, the
adjoint is uniquely determined up to natural equivalence. The basic
continuity properties also hold: left adjoints preserve colimits, and
right adjoints preserve limits. For ordinary categories, Freyd's
adjoint functor theorem is a partial converse to this result. It
says roughly\footnote{``Roughly'' because we are omitting all the size
  conditions in the statement: $\D$ should be locally small and when
  we say $G$ preserves limits we mean small limits.} that for $G : \D
\to \C$ to have a left adjoint it is sufficient that $\D$ is complete,
$G$ preserves limits and satisfies a further condition called the
\emph{solution set condition}. The precise form of the solution set
condition will not matter for us, only that (1) it is a size condition
in the sense that it requires there to exist a small set of morphisms
that somehow control all morphisms of a certain form, (2) it is
adapted to $G$: it is not just a condition on the category $\C$.

There is an adjoint functor theorem for \io-categories due to Lurie
which in some sense is less precise than Freyd's theorem for
$1$-categories in that its size condition is not adapted to $G$, but
rather is a global condition on $\C$ and $\D$. In practice, this is
not a problem: the conditions of Lurie's theorem are usually met when
they need to be.

\begin{theorem}[{\cite[Corollary 5.5.2.9]{HTT}}]
  Let $F : \C \to \D$ be a functor between \emph{presentable}
  \io-categories.
  \begin{enumerate}
    \item The functor $F$ has a right adjoint if and only if it
      preserves small colimits.
    \item The functor $F$ has a left adjoint if and only if it
      is \emph{accessible} and preserves small limits.
  \end{enumerate}
\end{theorem}

The terms ``presentable'' and ``accessible'' are what take the place
of the solution set condition (and other size conditions) in Freyd's
theorem. An \io-category is presentable if it has all small colimits
and is accessible. Accessibility is really the size condition; for
\io-categories it intuitively means that the \io-category is
determined by a small subcategory of objects which themselves are
compact, for functors between such categories it means it preserves
certain colimits which are a generalization of filtered colimits. For
precise definitions see \cite[Chapter 5]{HTT} or, for the analogous
theory in the case of ordinary categories, \cite{AdamekRosicky}, which
is highly recommended if only because it explains the relation of
these notions to universal algebra making them seem much less like
merely annoying technical set theoretic issues.

In classical category theory whenever $F : \C \to \D$ and $G: \D \to
\C$ are adjoint functors (with $F$ being the left adjoint), the
composite $G \circ F$ is a monad and $F \circ G$ is comonad. This is
also true in the \io-categorical context, but much harder to show
since, as the reader expects by now, the concept of monad requires not
an associative multiplication but one that is associative up to
coherent homotopy. The monad corresponding to an adjunction is
constructed in \cite[Section 6.2.2]{HA}; Lurie uses it to prove an
\io-analogue of the Barr-Beck theorem characterizing those adjunctions
in which $\D$ is equivalent to the category of $G\circ F$-algebras in
$\C$ and $G$ is equivalent to the forgetful functor. This result is
very useful in the theory of descent in derived algebraic geometry,
just as the classical version is useful in algebraic geometry.

Another construction, providing more explicit information, of the
monad and comonad corresponding to an adjunction will appear in
\cite{Adjunctions}. There the authors construct ``the walking
\io-adjunction'': an $(\oo,2)$-category $\A$ that has two objects $0$
and $1$, and two morphisms $f : 0 \to 1$ and $g:1 \to 0$ which are
adjoint to each other and form the free adjunction in the sense that
any pair of adjoint functors $(F,G)$ between \io-categories arise as
the images of $f$ and $g$ under some functor from $\A$ into the
$(\oo,2)$-category of all \io-categories. Given such a functor $H : \A
\to \Catio$ (where temporarily, $\Catio$ is an $(\oo,2)$-category),
the restriction of $H$ to $\Hom_\A(0,0)$ gives the monad $G \circ F$
with its multiplication and all the higher coherence data. This
$(\oo,2)$-category $\A$ is surprisingly just a $2$-category, the
``walking (ordinary) adjunction'' described in \cite{FreeAdj}.

\begin{remark}
  We mentioned above that if $F$ has a right adjoint $G$, then $G$ is
  canonically determined by $F$. In fact, all of the adjunction data
  is determined by $F$ in a sense similar to that in proposition
  \ref{comp}: there is an \io-category of adjunction data that has a
  forgetful functor to the arrow category of $\Catio$ which just keeps
  the left adjoint $F$; this forgetful functor has contractible
  fibers.  Similarly there is a contractible space of adjunction data
  with a given left adjoint $F$, right adjoint $G$ and unit $u$. The
  results in \cite{Adjunctions} describe more generally which pieces
  of adjunction data determine the rest up to a contractible space of
  choices.
\end{remark}

\subsection{Less basic \io-category theory}

Much more than just the basic notions of category theory have been
extended to \io-categories. There is a large number of topics we could
have included here, and we picked only two that are important in
applications: the different sorts of fibrations between
\io-categories, and stable \io-categories. The various fibrations play
a larger role in actual use of \io-categories, specially when
incarnated as quasi-categories, than do the notions in ordinary
category theory they generalize. Stable \io-categories are widely used
in derived algebraic geometry, since they replace triangulated
categories. Stable \io-categories are also downright pleasant to work
with and we have no compunctions about advertising them.

A very important topic we've omitted is the study of monoidal and
symmetric monoidal \io-categories. As usual in higher category theory
they're tricky to define since one must make the tensor product
associative up to coherent homotopy (or associative and commutative up
to coherent homotopy in the symmetric case). We've decided to omit
them since we feel that the coherence issues we'd point out about them
would be similar to the ones we've mentioned already in other
contexts. Also, there are gentler introductions to them than the one
in \cite{HA} already available: we recommend Lurie's old version in
\cite{DAG2} and \cite{DAG3}, or Groth's nice exposition \cite[Sections
3 and 4]{Groth}.

\begin{remark}
  On the topic of the point of view on monoidal categories in
  \cite{HA}, namely, that they are special cases of colored operads
  (note that \io-categories are also special cases), we mention the
  work of Ieke Moerdijk and his collaborators. They use dendroidal
  sets to model \io-colored operads instead of Lurie's simplicial
  sets. Just like simplicial sets are well adapted to taking nerves of
  categories, dendroidal sets are shaped to produce nerves of operads
  more naturally. They have developed dendroidal analogues of several
  of the models for $(\oo,1)$-categories we described and shown their
  equivalence, see \cite{Ieke} for a more leisurely description of this
  work than that available in the original sources.
  An upcoming paper will prove the equivalence between the dendroidal
  approach and Lurie's simplicial model \cite{Gijs}.
\end{remark}

Finally, a topic that will surely become increasingly important is the
theory of enriched \io-categories. Currently this is dealt with in a
somewhat ad hoc manner when it arises, but David Gepner and Rune
Haugseng are in the process of producing a systematic treatment
\cite{Rune}.

\subsubsection{Fibrations and the Grothendieck construction}

The \emph{Grothendieck construction} takes as input a weak $2$-functor
$F : \C^{\op} \to \Cat$ where $\C$ is an ordinary category (thought of
as a $2$-category with only identity $2$-morphisms), and $\Cat$  is
the $2$-category of categories. Even though $\Cat$ is a strict
$2$-category, it makes sense to consider a weak functor $F$, that is,
one that does not preserve composition on the nose, but rather comes
equipped with a natural isomorphism $F(g) \circ F(f)
\xrightarrow{\simeq} F(f \circ g)$ satisfying a coherence condition
for compositions of three morphisms. It produces as output a category
$\E$ and a functor $P : \E \to \C$. The functor $P$ produced is always
what is called a \emph{Grothendieck fibration}, and the Grothendieck
construction provides an equivalence of $2$-categories between
$\Fun(\C^{\op}, \Cat)$ and the $2$-category of Grothendieck fibrations
over $\C$. We won't give precise definitions here, but we'll
illustrate in an example.

\begin{example}
 Let $\C$ be the category $\Top$ of topological spaces and $F$ the
 functor that assigns to each space $X$ the category of vector bundles
 on $X$ (the morphisms from a bundle $E_1 \to X$ to a bundle $E_2 \to
 X$ are maps $E_1 \to E_2$ which form a commuting triangle with the
 projections to $X$ and which are linear on each fiber). For a
 continuous function $f : X \to Y$, $F(f)$ is given by pullback of
 vector bundles. For this functor, the Grothendieck construction
 produces a category $\E$ whose objects are vector bundles $E \to X$
 on arbitrary base spaces $X$ and whose morphisms from a bundle $E_1
 \to X_1$ to a bundle $E_2 \to X_2$ are commuting squares
 \centerline{
   \xymatrix @C=5mm @R=5mm {
     E_1 \ar[d] \ar[r]^e & E_2 \ar[d] \\
     X_1 \ar[r]^x & X_2 \\
   }
 }
 such that the map $e$ is linear on each fiber. The projection $P : \E
 \to \C$ simply forgets the bundles and keeps the bases. Now let's
 spell out what being a Grothendieck fibration means in this example.
 There are certain morphisms in $\E$ that are distinguished: the
 \emph{Cartesian} morphisms for which the above square is a pullback.
 Given any map $x : X_1 \to X_2$ in $\Top$ and any vector
 bundle $E_2$ on $X_2$, there is always a Cartesian morphism in $\E$
 with codomain $E_2 \to X_2$ that $P$ sends to $x$: namely, the one in
 which $E_1$ is the pullback $x^*(E_2)$. Every morphism $e'$ in $\E$
 such that $P(e')$ factors through $x$ can be factored through the
 morphism $e : x^*(E_2) \to E_2$.
\end{example}

The Grothendieck construction is very versatile. First, it clearly
allows one to deal with weak $2$-functors to $\Cat$ \emph{without}
leaving the world of $1$-categories. It is used this way in the theory
of stacks, for example. There one wants to replace the notion of
a sheaf of sets with a $2$-categorical version that has values in
categories, or more typically, just groupoids. Such a thing can be
defined in terms of weak $2$-functors to groupoids, but can also be
handled as Grothendieck fibration, which is the approach typically
taken. See for example Angelo Vistoli's notes \cite{Vistoli}. (By the
way, that also serves as a references for the details on the
Grothendieck construction we skipped above.)

But the Grothendieck construction also allows one to
calculate limits and colimits of functors to $\Cat$. There is a dual
version that given a covariant $F : \C \to \Cat$ produces what is called a
Grothendieck opfibration $\E \to \C$, and it's not too hard to show that:

\begin{enumerate}
  \item $\lim F$ is given by the category of \emph{coCartesian}
    sections of $P : \E \to \C$, that is, sections $\sigma : \C \to
    \E$ such that $P \circ \sigma = \id_C$ and $\sigma(f)$ is a
    coCartesian morphism in $\E$ for every morphism $f$ of $\C$.
  \item $\colim F$ is given by the localization
    $\E[\text{coCart}^{-1}]$ of $\E$ obtained by inverting all
    coCartesian morphisms.
\end{enumerate}

Here $\lim$ and $\colim$ denote what are sometimes called
pseudo-limits and pseudo-colimits which are the closest analogues in
the world of $2$-categories to homotopy limits and colimits. Indeed,
the reader who knows how to perform the Grothendieck construction
should compare it with the description of homotopy limits and colimits
in section \ref{holim}.

All of this is generalized to \io-categories in \cite[Chapter 2]{HTT}.
(See also Moritz Groth's notes \cite[Section 3]{Groth}).  Lurie
defines \emph{(co)Cartesian fibrations} corresponding to Grothendieck
(op)fibrations, and also discusses Joyal's notions of left and right
fibrations. As in the classical case, functors $\C^{\op} \to \Catio$
for an \io-category $\C$ are classified by Cartesian fibrations $\E
\to \C$; and right fibrations are the subclass of Cartesian fibrations
classifying functors that land in the subcategory of $\Catio$ of
$\oo$-groupoids.

The constructions described above for limits and colimits of
categories also generalize to \io-categories. See \cite[Proposition
3.3.3.1]{HTT} for limits of \io-categories; Corollary 3.3.3.4
specializes the previous result to left fibrations and thus provides a
construction for homotopy limits of spaces. On the colimit side, there
are corollaries 3.3.4.3 and 3.3.4.6 which construct colimits in
$\Catio$ and homotopy colimits of spaces respectively.

\subsubsection{Stable \io-categories}

Stable \io-categories are a wonderfully practical replacement for the
notion of triangulated category that fixes many of the problems in
their theory. Many prominent examples of triangulated categories, are
given to us almost by definition as homotopy categories of naturally
occurring \io-categories with nice properties that make them
homotopical analogues of abelian categories. These \io-categories are
the \emph{stable} ones we will shortly define. For example, the
derived category of an abelian category $\A$ is basically defined as
the homotopy category of the relative category of chain complexes in
$\A$ with weak equivalences given by the quasi-isomorphisms. Let's now
give the definition.

\begin{definition}
  An \io-category is \emph{pointed} if it has an object which is both
  initial and terminal. Such an object is called a \emph{zero object}
  and denoted $0$.
\end{definition}

Between any two objects $X$ and $Y$ in a pointed \io-category there is
a unique homotopy class of \emph{zero morphisms}, those that factor
through the zero object.

\begin{definition}
  A \emph{triangle} in a  \io-category is a commuting square of the
  form:

  \centerline{
    \xymatrix @R=5mm @C=5mm {
      X \ar[r]^{f} \ar[d] & Y \ar[d]^{g} \\
      0 \ar[r] & Z \\
    }
  }
  That is, a triangle is a homotopy between $g \circ f$ and the zero
  morphism.  If the square is a pushout square, $g$ is called the
  \emph{cofiber} of $f$; if the square is a pullback, $f$ is called
  the \emph{fiber} of $g$.
\end{definition}

Fibers and cofibers of maps play a role similar to kernels and
cokernels in abelian categories and to cocones and cones in
triangulated categories.  In terms of these concepts the definition of
stable \io-category is easy to give:

\begin{definition}
  A \emph{stable} \io-category is a pointed \io-category where every
  morphism has a fiber and a cofiber and every triangle is a pushout
  if and only if it is a pullback.
\end{definition}

That's it, that's the whole definition. Of course, the notion of
\io-category is much more complicated than the notion of $1$-category,
but after making the initial investment in \io-categories, the
definition of stable \io-category is simpler not only than the
definition of triangulated category but even than that of abelian
category.

For readers familiar with the definition of triangulated category, we
now explain how the homotopy category of any stable \io-category is
canonically triangulated. The triangles are defined to be to diagrams
of the form

\centerline{
  \xymatrix @R=5mm @C=5mm {
    X \ar[r]\ar[d] & Y \ar[r]\ar[d] & 0 \ar[d] \\
    0 \ar[r] & Z \ar[r] & U \\
  }
}

\noindent where both squares are pushouts. This makes the outer
rectangle a pushout too, and since for spaces homotopy pushouts
of that form produce suspensions, we call $U$ the suspension of $X$
and write $U \simeq \Sigma X$. Suspension is the translation
functor of the triangulated structure.

The advantage of stable \io-categories over triangulated begins to be
visible even from here: in a triangulated category the cone of a
morphism is determined up to isomorphism but \emph{not} up to a
canonical isomorphism, and this is because the universal property we
should be asking of the cone is homotopical in nature. In a stable
\io-category the cofiber is determined, as all colimits are, up to a
contractible space of choices which is exactly canonical enough in the
\io-world to make the constructions functorial.

As further evidence of the simplicity of stable \io-categories, notice
that proving the octahedral axiom in the homotopy category becomes
simply a matter of successively forming pushout squares, and putting
down a suspension whenever we see a rectangle with zeros in the bottom
left and top right corners:

\centerline{
  \xymatrix @R=5mm @C=5mm {
    X \ar[r]\ar[d] & Y \ar[r]\ar[d] & Z \ar[r]\ar[d] & 0 \ar[d] & \\
    0 \ar[r] & Y/X \ar[r]\ar[d] &  Z/X \ar[r]\ar[d] &
    \Sigma X \ar[r]\ar[d] & 0 \ar[d] \\
    & 0 \ar[r] & Z/Y \ar[r] & \Sigma Y \ar[r] &  \Sigma(Y/X) \\
  }
}

(Here we've used notation like $X \to Y \to Y/X \to \Sigma X$ for
the objects in a distinguished triangle.)

\begin{further}
  As usual, we've only scratched the surface of the theory. We justed
  wanted to advertise the practicality of stable \io-categories over
  triangulated categories and refer the reader to \cite[Chapter 1]{HA}
  for the development of the theory, or to Groth's presentation
  \cite[Section 5]{Groth}.
\end{further}

\section{Some applications}

This final section will briefly describe some uses of \io-categories
outside of higher category theory itself. (The characterization of
$(\infty,n)$-categories mentioned in section \ref{comparison} is an
excellent application of \io-categories within higher category
theory.) Each of these applications is a whole subject in itself and
we cannot hope to do any of them justice here, rather, we hope merely
to whet the readers' appetite and suggest further reading.
Unfortunately many great topics had to be left out, such as algebraic
$K$-theory, where the use of \io-categories has finally allowed to
produce a universal property of higher $K$-theory similar to the
familiar one satisfied by $K^0$: approximately, that $K$-theory is the
universal invariant satisfying Waldhausen's additivity theorem.  For
precise descriptions and proofs of this universal property see
\cite{BarwickK} and \cite{BlumbergK}.  Another proof of Waldhausen
additivity in the \io-setting can be found in \cite{FioreLueck} (but
without the universal property), and an \io-version of Waldhausen's
approximation theorem is proved in \cite{Fiore}.

\subsection{Derived Algebraic Geometry}

One motivation behind studying stacks in algebraic geometry is a
desire to ``fix'' quotients of schemes by group actions. For example,
if an algebraic group $G$ acts freely on a variety $X$, sheaves on the
quotient $X/G$ are the same as $G$-equivariant sheaves on $X$. But if
the action of $G$ is not free, then the relation between $G$-equivariant
geometry on $X$ and geometry on the quotient $X/G$ is complicated,
aside from the difficulty in saying precisely what the quotient means
in this case. Passing to stacks ``fixes'' the quotient, so that sheaves
on the quotient stack $[X/G]$ are again $G$-equivariant sheaves on $X$.

Just as stacks can be seen as a way to improve certain colimits of
schemes so that they always behave as they do in the nice cases,
\emph{derived algebraic geometry} can be regarded as fixing a dual
problem: correcting certain limits so that they behave more often
as they do in good cases. Namely, one motivation for derived
algebraic geometry is to ``improve'' intersections so that they
behave more like transverse intersections. Some intersections are
already improved by passing from varieties to schemes, for example
consider the case of two curves $X$ and $Y$ of degrees $m$ and $n$
in the complex projective plane. If they intersect transversely,
according to B\'{e}zout's theorem they always have $mn$ points of
intersection. When the intersections are not transverse there might
be fewer than $mn$ points of intersection, but as long as $X$ and
$Y$ have no common irreducible component, there is a natural
way to assign multiplicities to the points so that counted with
multiplicity there are still $mn$ points of intersection.

We can write the result in the transverse case in terms of cohomology
classes as $[X] \cup [Y] = [X \cap Y]$, where $[X], [Y] \in
H^2(\CC\P^2,\Z)$ are the fundamental classes of $X$ and $Y$, and
$[X \cap Y] \in H^4(\CC\P^2,\Z) = \Z$ is the fundamental class of the
$0$-dimensional variety consisting of the finitely many points of
intersection of $X$ and $Y$. The point of this example is that when
the intersection of $X$ and $Y$ is not transverse (but still assuming
$X$ and $Y$ do not share a common component), the same formula is true
provided we think of $X \cap Y$ as a scheme instead of a variety (and
define fundamental classes in an appropriate way). The dimension over
$\CC$ of the local ring of $X \cap Y$ at the points of intersection
gives the relevant multiplicities.

By passing from varieties to schemes we reach a context where
computing $[X] \cup [Y]$ does not require full knowledge of the pair
$(X,Y)$, just knowledge of the intersection $X \cap Y$, and the
intersection is again a scheme, the same kind of object as $X$ and $Y$.
So schemes ``fix'' the intersection of curves that don't share a
component, but they don't fix intersections of curves that do share
a component nor do they help intersections in higher dimensions.
If $X$ and $Y$ are now subvarieties of $\CC\P^n$ of complementary
dimension that intersect in a finite set of points, the dimension
of the local ring of $X \cap Y$ at a point of intersection is no
longer the correct multiplicity to make $[X] \cup [Y] = [X \cap Y]$
true; now the multiplicity is given by Serre's formula:
\[ (\text{multiplicity at } p\in X \cap Y) =
   \sum_{i\ge 0} \dim_\CC \Tor_i^{\O_{\CC\P^n,p}}(\O_{X,p},\O_{Y,p}). \]

Notice that this sum \emph{cannot} be computed from just the knowledge
of $X \cap Y$ as a scheme, since that only determines the $\Tor_0$ term,
but the same formula suggests a fix: we'd be able to compute the correct
multiplicities just from $X \cap Y$ if we could arrange for
\[\O_{X \cap Y, p} = \O_{X,p} \stackrel{\L}\otimes_{\O_{\CC\P^n,p}}
 \O_{Y,p}, \]
where $\otimes^\L$ denotes the derived tensor product, a chain complex
whose homology gives the relevant $\Tor$ groups. We can attempt, then,
to define \emph{derived schemes} by replacing commutative rings in the
definition of schemes by the chain complex version of commutative
algebras: commutative differential graded algebras or cdgas for short,
and to compute structure sheaves of intersections using derived tensor
products. This suggests jumping into \io-waters, since the derived
tensor product is only naturally defined up to quasi-isomorphism while
sheaves expect to have values determined up to isomorphism.

\begin{example}
  The higher categorical point of view also suggests itself by
  considering the other case we mentioned that schemes don't fix:
  intersecting two curves with a common component. For the simplest
  possible example let's intersect the $y$-axis in the affine plane
  $\CC^2$ with itself: the ring of functions we get is $\CC[x,y]/(x)
  \otimes \CC[x,y]/(x) = \CC[x,y]/(x,x)$. If, as for ordinary rings we
  had $\CC[x,y]/(x,x) = \CC[x,y]/(x)$ the intersection would be
  $1$-dimensional and we'd have no hope of a uniform intersection
  theory. We want instead to have that modding out by $x$ twice is
  different than doing it once, and makes the quotient ring
  $0$-dimensional again. The philosophy of higher categories suggests
  that instead of forcing $x$ to be $0$, we make it isomorphic to $0$;
  that way if we mod out twice, we can add different isomorphisms each
  time. Let's do that using cdgas for simplicity.  Let's start out
  with $\CC[x,y]$ with $x$ and $y$ in degree zero. To compute
  $A:=\CC[x,y]/\!/(x)$ (where we've used a double slash to avoid
  confusion with the ordinary quotient of rings) we add a generator
  $u$ in degree one with $du = x$, to make $x$ isomorphic (maybe
  ``homologous'' is a better term for isomorphism in this
  commutative $\oo$-groupoid $A$) to $0$. We can readily compute
  $H_0(A) = \CC[x,y]/(x)$ (ordinary quotient), and $H_i(A) = 0$ for
  $i\neq 0$. Now, we quotient by $x$ again, to get
  $B:=A/\!/(x)=\langle x_0,y_0,u_1,v_1 : du=dv=x \rangle$ (where the
  subindices on the generators indicate the degree). The non-zero
  homology is $H_i(B) = \CC[x,y]/(x)$ for $i=0,1$. This computes the
  correct $\Tor$ terms from Serre's formula:
  $\Tor^{\CC[x,y]}_i(\CC[x,y]/(x), \CC[x,y]/(x))$ is also given by
  $\CC[x,y]/(x)$ for $i=0,1$ and is $0$ otherwise. For similar but
  less trivial examples see the introduction of \cite{DAG5}.
\end{example}

\subsubsection{The cotangent complex}

Another motivation for derived algebraic geometry comes from
deformation theory. The cotangent complex $\L_X$ of a scheme $X$ over
a field controls its deformation theory through various $\Ext$ groups
of $\L_X$. Roughly speaking, it can be used to count how many
isomorphism classes of first order deformations $X$ has, and to count
how many non-isomorphic ways each first order deformation extends
to a second order deformation, and so on. But it seems initially to be
a purely algebraic creature. Grothendieck asked in 1968 for a
geometrical interpretation of the cotangent complex, and derived
algebraic geometry provides one. We'll only sketch how the groups
$\Ext^i(\L_{X,k}, k)$ are reified through the derived affine schemes.

For $i=0$, we'll only need ordinary rings, in fact. Let $D_0 = \Spec
k[\epsilon]/(\epsilon^2)$ be the walking tangent vector: it is a tiny
scheme whose underlying topological space has only one point, but that
additionally carries a direction so that morphisms $D_0 \to Y$
correspond to tangent vectors on $Y$.  We have that $\Ext^0(\L_{X,x},
k) \simeq \Hom_\ast(D_0, (X,x))$ (where $\Hom_\ast$ means base-point
preserving morphisms), but to obtain similar representations of the
higher $\Ext$, we need to pass to derived rings.

One of the great advantages of admitting nilpotent elements in the
rings used for algebraic geometry is that it allows for schemes like
$D_0$. A good intuitive picture to have for derived schemes is that
they are to ordinary schemes as those are to reduced schemes, that is,
derived schemes can have something like ``higher nilpotent elements''
in their structure sheaves. We are about to see an example of that
now. Let $D_i = \Spec k[\epsilon_i]/(\epsilon_i^2)$ where the ring is
now a cdga with generator $\epsilon_i$ in degree $i$ (thinking of the
chain complex as an abelian higher groupoid, it has a single
non-identity $i$-morphism). As is nicely explained in \cite{VezzosiL},
it turns out that $\Ext^i(\L_{X,x}, k) \simeq \pi_0 Map_\ast(D_i,
(X,x))$, where $\Map_\ast$ denotes a mapping space in an \io-category
of derived pointed schemes.

\subsubsection{Rough sketches of the definitions}

Now we can say approximately what derived schemes are: they are what
you get by taking some definition of scheme in terms of commutative
rings and performing two replacements: (1) replacing notions from
ordinary category theory appearing in the definition with the
corresponding definitions in \io-category theory, and (2) replacing
the category of commutative rings with an \io-category of generalized
commutative rings (that allow extraction of higher $\Tor$ groups from
some tensor product operation defined for them).

\begin{remark}
  Since passing to derived schemes and passing to stacks are meant
  to solve independent problems in the category of schemes, it is
  entirely possible and even desirable to generalize schemes in both
  directions at once to obtain the theory of derived stacks.
  The references we will mention actually treat derived stacks as well.
  Derived stacks are also higher stacks in the sense that
  the groupoids appearing in the theory of ordinary stacks are
  replaced by $\oo$-groupoids, but it also makes sense to consider
  \emph{underived} higher stacks as in \cite{SimpsonStacks}.
\end{remark}

There are choices for both what definition of schemes to generalize
and what sort of rings to use. One can use the definition of schemes
as locally ringed spaces locally isomorphic to affine schemes or the
point of view of the functor of points. For rings one could use
commutative differential graded algebras, simplicial commutative rings
or, for application to homotopy theory, $E_\oo$-ring spectra. Also,
given any commutative ring $R$, we can consider $R$-algebras of each
of those three kinds.

Let's talk about the choice of rings first. The three notions we
listed are related by functors of \io-categories $\text{SCR}_R \to
\text{CDGA}_R \to E_\oo\text{-}R\text{-}\text{Alg}$. The situation is
very simple if $R$ is a $\mathbb{Q}$-algebra: the second functor is an
equivalence, the first functor is fully faithful and its image
consists of connective cdgas: those whose homology groups are
concentrated in non-negative degrees. In general, it gets a little
messy: neither functor nor their composite is fully faithful.
Comparing free algebras might be illuminating: the free cdga on one
generator $x$ of degree $0$, say, is just the polynomial ring $\Z[x]$
concentrated in degree $0$. On the other hand, the free
$E_\oo$-$\Z$-algebra $A$ on $x$ is quite different, since the
multiplication on $A$ is not strictly commutative but only commutative
up to coherent homotopy. This means for example, that permuting the
$n$ factors in the product $x^n$, doesn't quite fix $x^n$, rather it
produces automorphisms of $x^n$, and in particular we get a non-trivial
homomorphism $\Sigma_n \to \pi_1(A,x^n)$ from the symmetric group on
$n$ letters. But in $\Z[x]$ the multiplication \emph{is} strictly
commutative and this $\Sigma_n$-action is trivial. It seems fair to
say that simplicial commutative rings are sufficient for applications
of derived algebraic geometry to algebraic geometry itself, while
$E_\oo$-ring spectra are mainly for applications to homotopy theory.

Now we'll deal with the choice of definition of scheme. Unlike what
happens for the choice of rings, the choice of style of definition
does not lead to different notions of derived scheme.  Both the
locally ringed space point of view and the functor of points view as
well as the equivalence between them were described by Jacob Lurie
in his PhD thesis or, more recently, in \cite{DAG5}. The latter also
contains a very general definition of derived geometric objects that
includes derived schemes with any of the above mentioned choices
of rings, but also things like derived smooth manifolds studied by
David Spivak in his PhD thesis (reworked into \cite{SpivakMfld})
and derived analytic spaces (further studied in \cite{DAG9}).
The rest of the DAG series of papers, \cite{DAG7}, \cite{DAG8},
\cite{DAG9}, \cite{DAG10}, \cite{DAG11}, \cite{DAG12}, \cite{DAG13},
\cite{DAG14} contain a wealth of information about algebraic geometry
in the derived setting that we don't have any space to describe.

Bertrand To\"{e}n and Gabriele Vezzosi favor the functor of points
in their work on homotopical algebraic geometry \cite{HAG1},
\cite{HAG2}. Their approach is based on an idea of Deligne's: it is
possible to construct ordinary stacks, say, by a categorical procedure
that receives very little outside input. Namely, one starts with the
symmetric monoidal category of abelian groups: commutative monoids for
the tensor product there are commutative rings, so we can get our
hands on the category of affine schemes. Next, one needs a second
piece of information: a topology on the category of affine schemes.
Using this one can define purely categorically stacks for the
topology. Deligne observed that one could attempt to do this for other
symmetric monoidal categories and choices of topology. To\"{e}n and
Vezzosi develop an \io-version of this idea starting from a symmetric
monoidal \io-category \cite{AGMon}.\footnote{They had their work cut
  out for them, since at the time the best available model for
  symmetric monoidal \io-categories were monoidal model categories,
  which are very rigid: they have an actual monoidal structure,
  that is, one associative up to coherent \emph{isomorphism} in the
  ordinary category on which the model structure is defined.}

Both approaches require a notion of \io-version of sheaf and of
\io-topos.  A sheaf with values in $\C$ on an \io-topos $\X$ is just a
functor $\X^{\op} \to \C$ that sends colimits
  in $\X$ to limits in $\C$. Just as in ordinary topos theory, for
any topological space $X$ there is a topos $\Shv(X)$ of sheaves on $X$
(in the $1$-categorical case these are sheaves of sets, in the
\io-case they are sheaves of $\oo$-groupoids), which completely
determines the space $X$ if it satisfies the mild technical condition
of being sober.  Sheaves on a space $X$ are just defined to be sheaves
on the \io-topos $\Shv(X)$, but the theory does not require derived
schemes or stacks to have underlying topological spaces and is
developed for \io-toposes.

\begin{remark}
  The notion of \io-topos is due to Charles Rezk who described them as
  model categories (and called them model toposes) in
  \cite{RezkTopos}.  Lurie developed their theory, using
  quasi-categories, in his book \cite{HTT}. To\"{e}n and Vezzosi have
  also written about \io-toposes using simplicial categories
  \cite{HAG1} and Segal categories \cite{SegalTopos}. For the reader
  familiar with the notion of \emph{elementary topos} in ordinary
  category theory, we should point out that \io-toposes only
  generalize the notion of Grothendieck toposes.  A Grothendieck topos
  can be succinctly defined as a localization of a presheaf category,
  that is as a category $\E$ that admits a functor $F : \E
  \hookrightarrow \Fun(\C^{\op},\Set)$ for some $\C$ such that $F$ is
  fully faithful and has left adjoint which preserves finite limits.
  This definition generalizes to \io-toposes by replacing $\Set$ with
  $\oo$-groupoids.  Like their $1$-categorical counterparts,
  \io-toposes can be characterized by analogues of Giraud's axioms.
  There is also an \io-version of the notion of site and sheaf on a
  site, and \io-categories of sheaves on sites provide examples of
  \io-toposes; but unlike the case for ordinary toposes, these
  examples are not all the \io-toposes.
\end{remark}

\begin{further}
  For a through overview of higher and derived stacks see
  \cite{ToenHD}. An early example of derived moduli spaces, before all
  the machinery was in place, can be found in \cite{Kontsevich}. A
  more recent work on derived moduli spaces is \cite{HAGDAG}.
  More information about the deformation theory
  aspects of derived algebraic geometry, can be found in the
  already mentioned \cite{DAG10}, but we also recommend Lurie's
  earlier ICM address \cite{LurieICM}
  which gives an expository account. To see derived algebraic
  geometry in action, see, for example:
  \begin{enumerate}
    \item David Ben-Zvi and David Nadler's work on derived loop spaces
      in algebraic geometry, (for example, \cite{BN}), or
    \item their paper \cite{BFN}, joint with John Francis,
      about the \io-category $QC(X)$
      of quasi-coherent sheaves on a derived stack $X$, where they
      show that given derived stacks $X$ and $Y$ over $k$ satisfying
      an appropriate finiteness condition, any $k$-linear
      colimit preserving functor $F : QC(X) \to QC(Y)$ is given by a
      quasi-coherent sheaf $K \in QC(X \times_k Y)$ by means of an
      \emph{integral transform} $F(S) := (\pi_Y)_* (K \otimes \pi_X^*
      S)$.
  \end{enumerate}
\end{further}

\subsubsection{Topological modular forms}

It would take us too far afield to describe topological modular
forms in any kind of detail, but for readers interested in homotopy
theory we want to at least mention this application of derived
algebraic geometry.

There is a beautiful relation due to Quillen \cite{QuillenFGL}
between generalized cohomology theories which have an
analogue of the theory of Chern classes, properly called
complex-oriented cohomology theories, and formal groups laws. A
formal group law is\footnote{Actually the definition we give
is for \emph{one-dimensional} formal groups laws.  There are
also $n$-dimensional formal group laws that are defined in the
same way but with $x,y,z$ denoting $n$-tuples of variables and
$F(x,y)$ denoting an $n$-tuple of power series.} a power series
$F(x,y) = x+y+(\text{higher order terms})$, with coefficents
in some commutative ring $A$, that satisfies $F(x,F(y,z)) =
F(F(x,y),z)$. Given any Lie group or algebraic group, the power
series expansion of the multiplication at the origin gives a formal
group law\footnote{Of the same dimension as the group.}, but not
every such law arises from this construction. There is a notion of
isomorphism of formal group laws given by changes of coordinates
$x' = g(x)$, $y' = g(y)$ where $g$ is a power series. If $A$ is a
field of characteristic $0$, all one-dimensional formal group laws
are isomorphic to $F(x,y) = x+y$, and there is a classification due
independently to Lazard \cite{Lazard} and Dieudonn\'e
\cite{Dieudonne} of one-dimensional formal groups laws over
algebraically closed fields of characteristic $p$.

A complex-oriented cohmology theory has an associated one-dimensional
commutative group law $F$, that gives the formula for the first
Chern class of a tensor product of line bundles in terms of the
Chern classes of the individual line bundles, $c_1(L_1 \otimes L_2)
= F(c_1(L_1), c_1(L_2))$. For ordinary cohomology, the formal group
law is the additive one, $F(x,y)=x+y$.  For complex $K$-theory,
it is the multiplicative formal group law, $F(x,y) = x+y+xy =
(1+x)(1+y)-1$. Both of these come from expanding the product
of a one-dimensional commutative algebraic group, the additive
group $\mathbb{G}_a$, and the multiplicative group $\mathbb{G}_m$,
respectively. Over an algebraically closed field there is only one
other kind of connected commutative one-dimensional algebraic group:
elliptic curves, whose associated formal group laws correspond to
\emph{elliptic cohomology} theories. We can think of assembling
all of the elliptic cohomology theories into something like a sheaf
of cohomology theories on the moduli space of elliptic curves. If
we had such a thing we could imagine that taking global sections
of it would produces a new cohomology theory that samples from
all elliptic cohomology theories at once. That cohomology theory
is (almost\footnote{For one thing, we actually need a sheaf on
the Deligne-Mumford compactification of the moduli space.}) what
$\text{tmf}$, the cohomology theory of topological modular forms,
is supposed to be.

There are serious technical difficulties in constructing it, since
the notion of ``sheaf of cohomology theories'' isn't well-behaved
enough to use for this purpose. Instead one can attempt to represent
the cohomology theories by spectra, for which the \io-notion of sheaf
we mentioned earlier works well. Lifting this sheaf of cohomology
theories to a sheaf of spectra turns out to be made easier,
paradoxically, by making the problem harder and lifting instead
to a sheaf of $E_\oo$-ring spectra, making $\text{tmf}$ the global
sections of an object in derived algebraic geometry. This lifting
was first achieved (and shown to be essentially unique) by Mike
Hopkins, Haynes Miller and Paul Goerss using an obstruction theory of
$E_\oo$-ring spectra. A second construction due to Jacob Lurie uses
more of the theory of derived algebraic geometry: he defines a moduli
problem for (derived) elliptic curves that can only be stated in the
derived setting and then uses his general representability criterion
to show that the moduli problem is represented by a derived
Deligne-Mumford stack, whose global sections then give $\text{tmf}$.

\begin{further}
  For a great introduction to the topic see Paul Goerss's
  Bourbaki Seminar talk \cite{Goerss} and the references therein.
  (That paper is also a good introduction to derived schemes.)
  Also recommended is Lurie's \cite{LurieTMF} where he outlines the
  second construction mentioned above.
\end{further}


\subsection{The cobordism hypothesis}

Michael Atiyah \cite{Atiyah} proposed a mathematical definition of
\emph{topological quantum field theory} (henceforth TQFT) inspired by
Graeme Segal's axioms \cite{SegalCFT}\footnote{The reader might prefer
  the more recently published \cite{SegalCFT2} which is an expanded
  version of an unpublished paper that circulated widely for many
  years.} for \emph{conformal field theory}. The definition uses the
category $\Bord_{\langle n-1,n\rangle}$ whose objects are closed
oriented $(n-1)$-manifolds and whose morphisms $M \to N$ are
diffeomorphism classes of \emph{bordisms}\footnote{Perhaps
  confusingly, there is no difference between bordisms and
  cobordisms. We are following someone's terminology in each
  case\dots} from $M$ to $N$, that is, $n$-manifolds $W$ whose
boundary is identified with $\overline{M} \sqcup N$. Here
$\overline{M}$ denotes $M$ with the opposite orientation and the
diffeomorphisms we consider are those fixing the boundary.  This
category has symmetric monoidal structure given by taking disjoint
unions (both at the level of $(n-1)$-manifolds and $n$-dimensional
bordisms). The brunt of Atiyah's axioms is that a TQFT is a symmetric
monoidal functor $F : \Bord_{\langle n-1, n\rangle} \to \Vect$, where
$\Vect$ carries the monoidal structure given by the tensor product of
vector spaces.

This definition implies that one can calculate the value of $F$ on
some bordism $W$ by chopping it up into a composition of simpler
bordisms and composing the images under $F$ of those pieces. For
example, when $n=2$, we could chop any surface into a combination of
pairs of pants, cylinders and spherical caps. This suggests that such
a functor is determined by just a handful of its values. This approach
becomes more and more complicated in higher dimensions because the
pieces which are required to build up all bordisms increase in number
and complexity. One might hope that a simpler theory could result from
allowing more general ways to cut up bordisms: Atiyah's definition
allows cutting along codimension $1$ submanifolds, but if we allowed,
say, arbitrary codimensional cuts, we could triangulate every
manifold. The cobordism hypothesis doesn't make it quite that simple
to cut up manifolds, but we will allow pieces of all dimensions.

John Baez and James Dolan's \emph{cobordism hypothesis} \cite{BaezDolan}
is about such extended TQFT's and does indeed say that they are determined
by very few of their values, namely, that they are determined by where
they send single points! This hypothesis has now been proved in the $n=2$
case by Jacob Lurie and Mike Hopkins and in general by Lurie. A detailed
and very readable sketch of the proof is available as \cite{LurieCob}. It
might be a little unfair to say that the cobordism hypothesis fulfills the
promise of a simpler theory mentioned in the previous paragraph, since
there were complicated foundational issues involved in both the precise
statement in Lurie's paper and in its proof. We'll say a little more
about these issues after stating a version of the cobordism hypothesis.

To allow chopping up manifolds into pieces of arbitrary dimension, we
replace the domain category of a TQFT by a higher category
$\Bord^{\fr}_n$ whose $k$-morphisms are $n$-framed\footnote{An
  $n$-framing on a $k$-manifold $M$ is a trivialization of its
  stabilized tangent bundle, $TM \oplus
  \underline{\mathbb{R}}^{n-k}$.} $k$-manifolds with corners, and
whose $(k+1)$-morphisms are $n$-framed bordisms between such. (The
framing is used to remove some ambiguity coming from the
diffeomorphism group of $\mathbb{R}^{m}$, it is important to
understand and to know what happens without it, but we leave the
explanations to the references.) What we've said so far is a little
sloppy. Baez and Dolan's original formulation involved an
$n$-category, so the description of the morphisms we gave is only
valid for $k<n$; for $k=n$, we must again take diffeomorphism classes
of bordisms.

A $\C$-valued framed $n$-dimensional extended TQFT then is a symmetric
monoidal functor from $\Bord^{\fr}_n \to \C$; where $\C$ is an
arbitrary symmetric monoidal $n$-category. We generalize to arbitrary
symmetric monoidal $n$-categories partly because there is no canonical
substitute for $\Vect$ in this case; finding a $\C$ that extends
$\Vect$ to an $n$-category is an interesting problem of it own.
Fortunately, we don't need any particular symmetric monoidal
$n$-category for the statement:

\textbf{Baez and Dolan's cobordism hypothesis:} A framed
$n$-dimensional extended TQFT with values in $\C$ is completely
determined by it's value on the point. Moreover, the value on the
point is always a \emph{fully dualizable} object of $\C$ and
there is a bijection between isomorphism classes of such TQFTs
and isomorphism classes of fully dualizable objects in $\C$.

We won't give the precise definition of fully dualizable object, for
that see \cite[Section 2.3]{LurieCob}, but we will describe them
briefly. First, notice that the usual definition of adjoint functors
can be phrased in terms of $1$-morphisms and $2$-morphisms in $\Cat$
and thus makes sense in any $2$-category, and can even be used for
$k$- and $(k+1)$-morphisms in a higher category to define adjunctions
of $k$-morphisms. Call a $k$-morphism dualizable if it has both a left
and a right adjoint, and call it $n$-times dualizable if it is
dualizable and the $(k+1)$-morphisms that show it has adjoints are
themselves $(n-1)$-times dualizable\footnote{Where, of course,
  ``$1$-times dualizable'' just means ``dualizable''.}. Finally, a
monoidal $n$-category $\C$ can be thought of as an $(n+1)$-category
with a single object $B\C$, so we can say that an object of $\C$ is
$m$-times dualizable if as a $1$-morphism in $B\C$, it is $m$-times
dualizable. Then a fully dualizable object of a symmetric monoidal
$n$-category $\C$ is just an $(n-1)$-times dualizable object. Notice
that this really does depend on $n$: a symmetric monoidal $n$-category
can also be regarded as a symmetric monoidal $(n+1)$-category (with
only identity $(n+1)$-morphisms), but being $(n-1)$-times dualizable
is different from being $n$-times dualizable.

This $n$-category version of the cobordism hypothesis was proven for
$n=2$ by Chris Schommer-Pries in his PhD thesis \cite{S-PTQFT}, but
for higher $n$ a lack of a solid, practical theory of $n$-categories
impeded progress. Lurie's key insight was to prove instead a more general
version using the easier theory of $(\oo,n)$-categories and then
deduce the original formulation by truncation. So Lurie replaced the
$n$-category described above by one in which
\begin{enumerate}
  \item $n$-morphisms are bordisms, not diffeomorphisms classes of such,
  \item $(n+1)$-morphisms are diffeomorphisms,
  \item $(n+2)$-morphisms are isotopies of diffeomorphisms,
  \item $(n+3)$-morphisms are isotopies of isotopies, and so on.
\end{enumerate}

\begin{further}
  For a more detailed introduction to the cobordism hypothesis and its
  applications, Dan Freed's survey \cite{Freed} is highly recommended,
  as is, of course, Lurie's \cite{LurieCob} which is not just a
  detailed outline of the proof but also includes a lot of motivation
  (and an introduction to higher category theory!), and many other
  interesting versions of the homotopy hypothesis, and applications.
  The reader might also be interested in Julie Bergner's survey
  \cite{BergnerCob} which focuses on models for $(\oo,n)$-categories
  but does state Lurie's version of the cobordism hypothesis in some
  detail, describing the construction of $\Bord^{\fr}_n$.
\end{further}

\bibliographystyle{amsalpha} \bibliography{infinity-survey}{}
\end{document}